\documentclass[reprint,
 aps,
]{revtex4-2}

\newtheorem{proposition}{Proposition}
\newtheorem{theorem}{Theorem}
\newtheorem{lemma}{Lemma}
\usepackage{mathtools, amsmath, amssymb, graphicx, amsfonts}
\usepackage{xcolor} %to color comments
\usepackage{bbm} %to use indicator function
\usepackage[normalem]{ulem} %to strikeout text (\sout)
\usepackage{float} %to use [h!] with figures
\usepackage{comment} %to comment out parts of the text
\usepackage{amsopn} %standard?

\usepackage[breaklinks=true]{hyperref}
\usepackage{breakcites}

\begin{document}

\preprint{APS/123-QED}

\title{Non-linear network dynamics with consensus-dissensus bifurcation}

\author{Karel Devriendt}
 \email{devriendt@maths.ox.ac.uk}
 \altaffiliation[Also at ]{Alan Turing Institute, London, UK}
\author{Renaud Lambiotte}%
\affiliation{%
 Mathematical Institute, University of Oxford, Oxford UK }%

\date{\today}

%%%%%%%%%%
% ABSTRACT
%%%%%%%%%%
\begin{abstract}
We study a non-linear dynamical system on networks inspired by the pitchfork bifurcation normal form. The system has several interesting interpretations: as an interconnection of several pitchfork systems, a gradient dynamical system and the dominating behaviour of a general class of non-linear dynamical systems. The equilibrium behaviour of the system exhibits a global bifurcation with respect to the system parameter, with a transition from a single constant stationary state to a large range of possible stationary states. Our main result classifies the stability of (a subset of) these stationary states in terms of the \emph{effective resistances} of the underlying graph; this classification clearly discerns the influence of the specific topology in which the local pitchfork systems are interconnected. We further describe exact solutions for graphs with external equitable partitions and characterize the basins of attraction on tree graphs. Our technical analysis is supplemented by a study of the system on a number of prototypical networks: tree graphs, complete graphs and barbell graphs. We describe a number of qualitative properties of the dynamics on these networks, with promising modeling consequences.
\end{abstract}

\maketitle

%%%%%%
%BODY
%%%%%%
\section{Introduction}\label{S1: introduction}\noindent
Network dynamics are widely used as a natural way to model complex processes taking place in systems of interacting components. Within this framework,  time-varying states are assigned to the nodes of a network  and evolve according to  interaction rules defined between neighbouring nodes. Sufficiently simple  for theoretical investigations, the resulting dynamics may yet exhibit complex emergent behaviour of the global network state, making them suitable to model various real-world systems. Moreover, the interplay between the underlying network structure and the rich phenomenology of dynamics taking place on it make network dynamics a powerful tool to better understand and characterize the network itself. Some well-known examples of network dynamics include random walks \cite{Lovasz, pvm_PA}, epidemic spreading \cite{pvm_epidemics}, synchronization of oscillator systems \cite{Strogatz_review,Arenas}, consensus dynamics and voter models \cite{Saber,Castellano} and power grids \cite{Dorfler_electrical}. An overview of these applications and many other examples can be found in \cite{Barrat, gleeson_porter, Strogatz_chaos}.
\\
In this article, we propose a new non-linear dynamical system inspired by the pitchfork bifurcation normal form. Our choice of dynamical equations is supported by a number of different interpretations. We find that the system can be seen as (i) a set of interacting (1D) pitchfork systems, (ii) a gradient dynamical system for a potential composed of double-well potentials over the links of the network and finally (iii) as the dominating behaviour of a general class of non-linear dynamics with odd coupling functions. Qualitatively, the main property of the system is that it exhibits a \emph{bifurcation} in the possible stationary states. In the first parameter regime, our system is essentially diffusive and evolves to a unique, uniform stationary state. In the second parameter regime the coupling function is a mixed attractive/repulsive force and the equilibrium is characterized by a large number of stationary states. We find an explicit description for (a subset of) these stationary states and analyse their stability using linear stability analysis. Our main technical result  classifies the stability of these stationary states in terms of the \emph{effective resistance} of certain links. The effective resistance is a central concept in graph theory with links to random walks \cite{Lovasz, Doyle}, distance functions and embeddings \cite{Klein, Fiedler_book, krl_simplex}, spectral sparsification \cite{Spielman} and many more. Its appearance as a determinant for (in)stability in our non-linear dynamical system is very surprising and at the same time a perfect example of the rich interplay between structure and function in network dynamics. Furthermore, analytical results are found for the basins of attraction (on tree graphs) of the stationary states, and an exact solution of the system is derived for certain types of graphs which include graphs with external equitable partitions. The latter result adds to a long list of interesting observations of dynamics on graphs with (external) equitable partitions and related symmetries \cite{Schaub, Pecora, Bonaccorsi, krl_umff, Ashwin_swift, Golubitsky_groupoid_formalism}. 
\\
Our technical analysis is supplemented by a detailed description of the system on complete and barbell graphs. On the complete graph, we find that a subset of the stable stationary states determine a \emph{balanced} bipartition of the graph with each group corresponding to one of two existing state values, and neither group being too dominant (hence balanced). On the barbell graph, a similar balanced bipartition is observed within each of the complete components but with a non-zero difference between the average states of both components. We discuss how these observations might be interpreted in the framework of opinion dynamics.
\\~
\\
Our choice to focus on a specific dynamical system is restrictive in various ways and our results only pertain to a small corner of the theory of non-linear dynamics on networks as a consequence. \textcolor{black}{In a follow-up on the present work however, we found that our results generalize to a much broader class of non-linear systems \cite{Marc_arxiv}, suggesting a potential wider relevance.} Other works on this subject, notably the results of Golubitsky \emph{et al.} \cite{Golubitsky_groupoid_formalism, Golubitsky_Gandhi} and Nijholt \emph{et al.} \cite{Nijholt_thesis, Nijholt_center_manifolds} describe and characterize general classes of systems whose dynamics are constrained by a given underlying structure. Their results allow to determine which dynamical features (e.g synchronization conditions, bifurcations) are robust (generic) with respect to the network structure; in other words it details which features can be explained purely from the network structure irrespective of the specific choice of coupling functions. Our contributions are no attempt at such generality \textcolor{black}{on the system level}, but instead aim at developing a qualitative understanding of non-linear dynamics on graphs, starting from a basic ‘toy’ system and describing its interesting properties, with a focus on the influence of the network structure on these properties. \\
\textcolor{black}{A second relevant line of research is the recent work by Franci \emph{et al.} \cite{Franci} and Bizyaeva \emph{et al.} \cite{Bizyaeva} which study decision-making in (multi-option) opinion dynamics. They formulate opinion dynamics in a fully generalized setting, and show - independent of further system models - that this setting can exhibit a variety of rich non-linear dynamical features such as consensus-dissensus bifurcations and opinion cascades. The model analysed in our article fits in the framework of Franci \emph{et al.} as a particular two-option opinion dynamical system and consequently, certain features such as the global consensus-dissensus bifurcation and the observations in Section \ref{S6: examples} can be explained in this context. However, our contributions are complementary to those made in \cite{Franci, Bizyaeva} as our particular model choice allows us to derive many other specific and interesting results, in particular related to the stability of stationary states and exact solutions in the presence of external equitable partitions. Furthermore, the main results of Franci \emph{et al.} follow from a so-called equivariant analysis of the system which deduces properties of the system, starting form its symmetries. Our results (and those in \cite{Marc_arxiv}) follow from an algebraic and graph theoretic analysis instead and are valid for more general network structures as a result.}
\\
\textcolor{black}{A third} body of related work is the well-developed field of coupled oscillator systems \cite{Strogatz_review, Ashwin_swift, Arenas}, where many similar questions are studied for non-linear (oscillator) systems on networks. In Section \ref{S7: synchronization} we briefly discuss the setup of coupled oscillator systems and highlight a particular result from \cite{Dorfler_smartgrids} which closely relates to our stability result, Theorem \ref{th: single link stability}.
\\~
\\
The rest of this paper is organized as follows.
Our dynamical system is introduced in Section \ref{S2: non-linear system} together with a number of interpretations of the system. Section \ref{S3: stationary states} introduces the notion of stable and unstable stationary states, and describes the stability results for our system. Section \ref{S4: exact solutions} describes some cases where the system equations can be solved exactly, and Section \ref{S5: basins of attraction} deals with the characterization of basins of attraction. In Section \ref{S6: examples} finally, system \eqref{eq: dynamical system} is studied on a number of prototypical networks with a focus on the qualitative behaviour of the solutions. A related result about synchronization in coupled oscillators is described in Section \ref{S7: synchronization}, and the article is concluded in Section \ref{S8: conclusion} with a summary of the results and perspectives for future research.
%%%%
%
%
%%%%
\section{The non-linear system}\label{S2: non-linear system}
We will study a dynamical system defined by a set of non-linear differential equations that determine the evolution of a dynamical state $\mathbf{x}(t)$. This state is defined on a graph $G$ where each of the $N$ nodes has a corresponding state value $x_i(t)\in\mathbb{R}$ which together make up the system state as $\mathbf{x}(t)=(x_1(t),\dots,x_N(t))$. The dynamics of $\mathbf{x}(t)$ are determined at the node level by a non-linear \emph{coupling function} between neighbouring nodes. For a node $i$ with neighbours $j\sim i$, the dynamics are described by
\begin{equation}\label{eq: dynamical system}
\frac{dx_i}{dt} = \sum_{j\sim i} r(x_i-x_j) - (x_i-x_j)^3
\end{equation}
where $r$ is a scalar parameter, called the \emph{system parameter}. Since the states are coupled via their differences, the average state value does not affect the dynamics and the \emph{state space} of system \eqref{eq: dynamical system} is thus equal to $X=\mathbb{R}^N/\mathbf{1}$, i.e. with any two states $\mathbf{x}$ and $\mathbf{y}$ equivalent if $\mathbf{x}-\mathbf{y}$ is constant for all nodes. In other words, the dynamics is translation invariant. When considering a specified initial condition $\mathbf{x}(0)=\mathbf{x}_0$, we will also write the solution of system \eqref{eq: dynamical system} as $\mathbf{x}(t,\mathbf{x}_0)$.
\\
\textcolor{black}{There are various ways to interpret the node states. In the setting of (linear) consensus dynamics, as used frequently in the robotics and control community, the state variable $x_i(t)$ represents a real-valued parameter or measurement of an agent in a physical system and the goal is to coordinate these variables globally by following some local dynamics \cite{Saber}, similar to our system \eqref{eq: dynamical system}. In the setting of opinion dynamics \cite{Degroot, Franci} on the other hand, the node states $z_i(t)\in\mathcal{I}$ in an interval (usually $\mathcal{I}=[0,1]$) reflect the commitment of an agent in the network to an option/belief A ($z_i=0$) or to an alternative B ($z_i=1$) instead \footnote{\textcolor{black}{In the more general case of multiple, say $N_o$ options the state of every node $i$ corresponds to a point in the $N_o$ simplex $\mathbf{z}_i(t)\in\Delta_{N_o}$. Even more general, the celebrated Degroot model \cite{Degroot} features a distribution (measure) on the simplex as a node state.}}; the state dynamics then model the (social) processes by which agents update their opinions or beliefs. As shown in \cite{Bizyaeva}, there is a mapping of (forward invariant and bounded) dynamics in $\mathbb{R}^N$, i.e. system \eqref{eq: dynamical system} on $X$, to opinion dynamics with state space $Z$. In this context, the system parameter $r$ is sometimes interpreted as a measure of social attention or susceptibility to social influence. Our system can thus be seen as a non-linear generalization of consensus dynamics (see also \cite{Srivastava}) or can be mapped onto a two-option opinion dynamics model; the further derivations in this article will be independent of these interpretations.}
\\~\\
In what follows, we show how our system appears naturally in three different settings. Apart from suggesting different motivations for the study of our system, each perspective comes with a set of tools and results that will be used in our further analysis.
\subsection{Three perspectives on the dynamical system}\label{SS1A: three perspectives}
\subsubsection{Pitchfork bifurcation normal form}
\begin{figure}[h!]
    \centering
    \includegraphics[scale=0.45]{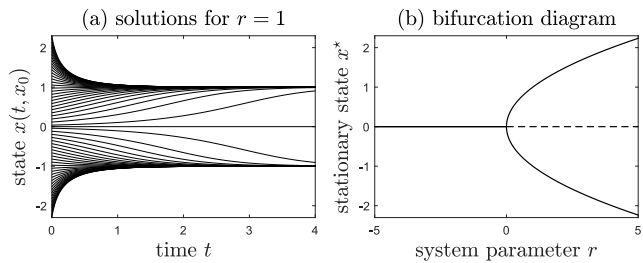}
    \caption{(a) Exact solutions for the pitchfork system $dx/dt=p(x)$ (as described in Section \ref{S4: exact solutions}, Appendix \ref{A3: exact solution}) for $r=1$ and a range of initial conditions $x_0$. These solutions illustrate the stable $(x^\star=\pm\sqrt{r})$ and unstable $(x^\star=0)$ stationary states for positive $r$. (b) Together with the stable stationary state $x^\star=0$ for negative $r$, these solutions determine the characteristic, eponymous bifurcation diagram of system \eqref{eq: dynamical system}}
    \label{fig: pitchfork}
\end{figure}

The definition of system \eqref{eq: dynamical system} is inspired by the so-called \emph{pitchfork bifurcation} dynamical system. This $1$-dimensional system with state $x(t)\in\mathbb{R}$ is given by the non-linear differential equation 
\begin{equation}\label{eq: 1d pitchfork system}
\frac{dx}{dt} = rx - x^3 \text{~with parameter~}r,
\end{equation}
where we will further also use the short-hand notation $p(x)=rx-x^3$ for the pitchfork function. System \eqref{eq: 1d pitchfork system} is the prototypical form (i.e. normal form) for dynamical systems that exhibit a bifurcation from a single stationary state to three distinct stationary states \cite{Strogatz_chaos}. This bifurcation occurs between a single stable stationary state $x^\star=0$ when $r<0$, and two stable states $x^\star=\pm\sqrt{r}$ and one unstable state $x^\star=0$ when $r>0$. Figure \ref{fig: pitchfork} shows the solutions of the pitchfork system (see also Section \ref{S4: exact solutions}) and illustrates the characteristic \emph{bifurcation diagram} to which the system thanks its name.

The system studied in this article thus consists of a pitchfork bifurcation system for the state difference $(x_i-x_j)$ over each of the links, with interactions coming from the shared variables of links with common nodes. Unsurprisingly, the larger interconnected system exhibits more complex behaviour than each of the smaller systems added together. In particular, our main result Theorem \ref{th: single link stability} highlights that the stable stationary states of the interconnected system can differ greatly depending on the way in which the links are interconnected. 
\\
Another way to see that system \eqref{eq: dynamical system} is closely related to the pitchfork bifurcation normal form is by introducing the \emph{link variable} $y_{{\ell}}=(x_i-x_j)$ for all links $\ell=(i,j)\in\mathcal{L}$ with an orientation $i\rightsquigarrow j$ fixed by taking the difference $(x_i-x_j)$. The dynamics can then be rewritten as
\begin{equation}\label{eq: system with link variables}
\frac{d y_{\ell}}{dt} = 2p(y_{\ell}) + \sum_{{m}\sim{\ell}}\sigma({\ell},{m})p(y_{m})
\end{equation}
where two links ${m}\sim{\ell}$ meet if they share a common node, and where the sign $\sigma({\ell},m)=\pm1$ of the interaction term depends on the relative orientation of the links; the matrix with entries $\sigma(\ell,m)$ for adjacent links and zero otherwise is also referred to as the edge adjacency matrix.
\subsubsection{Dominating behaviour of odd coupling functions}
System \eqref{eq: dynamical system} is a specific example of a more general class of non-linear dynamical systems on a graph:
\begin{equation}\label{eq: odd non-linear system}
\frac{dx_i}{dt} = \sum_{j\sim i}f(x_i-x_j) \text{~with odd function~}f.
\end{equation}
An important property of this class of systems is that the average state $\langle \mathbf{x}\rangle \triangleq \tfrac{1}{N}\sum x_i$ is always a \emph{conserved quantity} \footnote{The conserved quantity \unexpanded{$\langle \textcolor{black}{\mathbf{x}}\rangle$} originates from the symmetry (around the origin) of the coupling function, in close resemblance to Noether's celebrated connection between conservation laws and symmetries.} for the dynamics. If we furthermore assume $f$ to be analytic, the dominating behaviour for systems of the form \eqref{eq: odd non-linear system} around the consensus state can be studied by looking at the Taylor expansion of $f$ around $(x_i-x_j)=0$ as
$$
\frac{dx_i}{dt} = \sum_{j\sim i}\frac{df}{dx}\Big\vert_{0}(x_i-x_j) + \frac{1}{6}\frac{d^3f}{dx^3}\Big\vert_{0}(x_i-x_j)^3 + O((x_i-x_j)^5).
$$
A first-order approximation retrieves a simple, linear diffusion process. For the third-order approximation on the other hand, we see that by introducing the parameter $r=-6(\tfrac{df}{dx}\big/\tfrac{d^3f}{dx^3})\vert_{\mathbf{x}=0}$ and rescaling time as $t' = -(6\big/\tfrac{d^3f}{dx^3})t$ we retrieve system \eqref{eq: dynamical system}. In other words, the analysis of system \eqref{eq: dynamical system} is indicative for a general class of non-linear systems with odd coupling functions in the near-consensus regime \footnote{A given function $f$ fixes the value of $r$ which means the detailed-balance stationary states with state differences $\sqrt{r}$ might be far from consensus. For instance, for $f:x\mapsto \sin(x)$ we have $r=6$ for which the difference over dissensus links will be $(x^\star_i-x^\star_j)=\sqrt{r}\approx 2.45$, which is far from the consensus value $0$ and thus makes the approximation inaccurate.}.
\\
In \cite{Srivastava}, systems of the form \eqref{eq: odd non-linear system} are considered within the general problem of non-linear consensus and called relative non-linear flow. They are studied alongside absolute non-linear flow, of the form $dx_i/dt = \sum \left( f(x_i)-f(x_j)\right)$ and disagreement non-linear flow, of the form $dx_i/dt = f(\sum(x_i-x_j))$. While some general results are found for the latter two, the discussion of relative flow systems in \cite{Srivastava} is limited to the description of a number of small systems.
%%%
%
%%%
\subsubsection{Gradient dynamical system}
\begin{figure}[h!]
    \centering
    \includegraphics[scale=0.6]{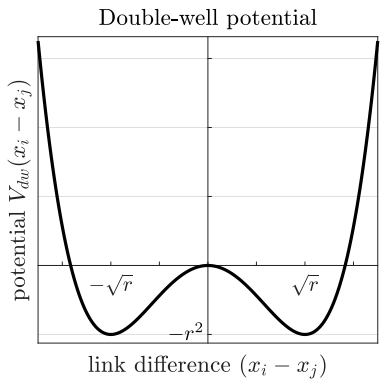}
    \caption{A double-well potential is a symmetric potential function with two local minima (the `wells') separated by a local maximum. In the case of our gradient system \eqref{eq: dynamical system}, the potential function $V$ is composed of a double-well potential term for each of the link differences, as $V(\mathbf{x})=\sum_{j\sim i}V_{\text{dw}}(x_i-x_j)$ where the specific double-well function $V_{\text{dw}}$ following from \eqref{eq: potential function} is illustrated above.}
    \label{fig: potential well}
\end{figure}

System \eqref{eq: dynamical system} also has the strong property that it is a \emph{gradient dynamical system}. This means that there exists a \emph{potential function} $V:X\rightarrow \mathbb{R}$ on the state space, such that the state dynamics are given by the negative gradient of this potential. For system \eqref{eq: dynamical system}, the potential takes the form 
\begin{equation}\label{eq: potential function}
V(\mathbf{x}) = \frac{1}{4}\sum_{i\sim j}(x_i-x_j)^2((x_i-x_j)^2-2r),
\end{equation}
from which the dynamics are retrieved as $d\mathbf{x}/dt = -\nabla V(\mathbf{x})$ with the gradient operator $\nabla=\left(\frac{\partial}{\partial x_1},\frac{\partial}{\partial x_2},\dots,\frac{\partial}{\partial x_N}\right)^T$. Interestingly, we see that the potential $V(\mathbf{x})$ in \eqref{eq: potential function} is composed of a separate potential term for each of the links. As illustrated in Figure \ref{fig: potential well} these terms are equal to a \emph{double-well potential}, which are minimal at $(x_i-x_j)=\pm \sqrt{r}$ separated by a local maximum at $x_i=x_j$. As we will see later, the link differences at these local optima also appear as stationary solutions of the system. 
\\
An important feature of gradient dynamical systems is that the potential is a decreasing function of time, i.e. the potential satisfies $\dot{V}\leq 0$ with equality if and only if the system is at a stationary point. This means that system \eqref{eq: dynamical system} is \emph{dissipative} for the potential, in contrast with its \emph{conservation} of the average value $\dot{\langle \mathbf{x}\rangle}=0$. This feature restricts the possible evolution of a system for a given initial state as $V(\mathbf{x}(t,x_0))\leq V(\mathbf{x}_0)$ must always be satisfied.
\section{Stationary states}\label{S3: stationary states}
Starting from the definition of the dynamics \eqref{eq: dynamical system}, we  study a number of different aspects of the system. A first important characterization is the long-term behaviour of the dynamics: starting from some initial state at $t=0$, in which states can we expect to observe the system after waiting sufficiently long? This question is answered by studying the \emph{stationary states} $\mathbf{x}^\star$ of the system, which are equilibrium states where the system is at rest, i.e. characterized by $d\mathbf{x}^\star/dt=0$. Since these states are the only points in which the potential does not strictly decrease, the system is guaranteed to evolve to a stationary state eventually. From a practical perspective, the stronger notion of \textcolor{black}{\emph{asymptotically stable}} stationary states is interesting. These are states for which the system is in a robust equilibrium, i.e. in the case of some perturbation $\boldsymbol{\epsilon}$, the state $\mathbf{x}^\star+\boldsymbol{\epsilon}$ will evolve back to $\mathbf{x}^\star$. We start by determining (a subset of) the stationary states of system \eqref{eq: dynamical system}, followed by an analysis of their stability.
\\
A direct translation of the stationarity condition yields the following characterization of stationary states:
\begin{equation}\label{eq: stationary states system}
\frac{d\mathbf{x}^\star}{dt}=0 \Leftrightarrow \sum_{j\sim i} r(x^\star_i - x^\star_j) - (x^\star_i - x^\star_j)^3=0,~\forall i.
\end{equation}
Generally, finding a stationary state $\mathbf{x}^\star$ thus involves solving a (potentially large) system of cubic equations. However, the possible solutions for $\mathbf{x}^\star$ differ greatly depending on the value of $r$. When $r<0$, only a single stationary state is possible: the \emph{consensus stationary state} where each node state equals the same constant value ${x}^\star_i = c$, equivalent to $\mathbf{0}\in X$ in the state space. In other words, for $r<0$ system \eqref{eq: dynamical system} is a (non-linear) form of diffusion or consensus dynamics. For $r>0$ on the other hand, the equations for $\mathbf{x}^\star$ can have many different solutions. Consider the case where a pair of linked nodes $i$ and $j$ have a state difference equal to $\sqrt{r}$. Then the coupling function between $i$ and $j$ will vanish, since $p(\sqrt{r})=0$, and the same happens when this difference is equal to $-\sqrt{r}$ or $0$. As a consequence, if the difference over all links equals one of these three values, all of the coupling functions will be inactive and the system will be in equilibrium. In other words, any state of the form
\begin{equation}\label{eq: detailed-balance states}
\mathbf{x}^\star \text{~with~} (x_i^\star - x_j^\star)\in\left\lbrace 0,\pm\sqrt{r}\right\rbrace \text{~for all~} i\sim j
\end{equation}
is stationary. Since the global equilibrium in these states originates from a local equilibrium (balance) for each of the links, we refer to solutions of this type as \emph{detailed-balance} stationary states \textcolor{black}{\footnote{\textcolor{black}{This terminology is inspired by chemical reaction kinetics, where detailed balance refers to the fact that a chemical system is in equilibrium (globally) when all individual microscopic processes, i.e. the occurring chemical reactions in the system, are balanced (locally). Similarly, in a detailed-balance stationary state $\mathbf{x}^\star$ we achieve a global equilibrium $d\mathbf{x}^\star/dt=0$, due to the conditions $p(x^\star_i-x^\star_j)=0$ for each of the links locally.}}}. Links with a zero difference will also be called \emph{consensus} links and links with a $\pm\sqrt{r}$ difference \emph{dissensus} links. \textcolor{black}{We note that the terms consensus and dissensus are often used to describe the global state of a system instead, while two neighbouring nodes having the same (different) state is then called (dis)agreement \cite{Saber, Franci}. Our choice to use the consensus/dissensus terminology at the level of single links follows from our interpretation of system \eqref{eq: dynamical system} as an interconnected collection of smaller systems, each of which can be in consensus, dissensus or another state (in non-detailed balance states). In certain cases these two uses of consensus and dissensus coincide, see e.g. Section \ref{SS6B: complete graph}.}
\\~\\
We recall from Section \ref{SS1A: three perspectives} (and Figure \ref{fig: potential well}) that a dissensus link corresponds to a minimum for the double-well potential and a consensus link to a local maximum. This means that local stationary states are composed to form global stationary states. In principle, $3^L$ possible detailed-balance solutions exist, with each link independently taking one of the possible differences. When the graph contains cycles however, these differences must be consistent across each cycle which reduces the number of possible solutions, down to a minimum of just $2^N$ possible detailed-balance states for the (maximally cyclic) complete graph (see Section \ref{SS6B: complete graph}).\\
From the perspective of gradient dynamics, the potential of detailed-balance stationary states can be expressed compactly in terms of the \emph{number of dissensus links} $\ell$ as
\begin{equation}\label{eq: potential in terms of dissensus links}
V(\mathbf{x}^\star) = -\tfrac{1}{4}r^2\ell.
\end{equation}
In other words, the higher the fraction of dissensus links in a stationary state, the lower the corresponding potential. We use this result in Section \ref{S5: basins of attraction} when describing the basins of attraction on tree graphs.
\\~\\
When solving equations \eqref{eq: stationary states system} directly or simulating the system, other stationary states can be found. \textcolor{black}{In the case of highly symmetrical graphs for instance, tools from equivariant dynamics can be used to find explicit descriptions of stationary states \cite{Golubitsky_advances_in_symmetry}. Generally,} whenever a graph has cycles (i.e. it is not a tree as in Section \ref{SS6A: loopless networks}) solutions may exist which are not detailed-balance states. Such states are difficult to describe in general and might even be degenerate. On the $3$-cycle graph for instance, all states $\mathbf{x}^\star$ on the circle
$$
\left\lbrace \mathbf{x} \in X : (x_1-\langle \mathbf{x}\rangle)^2+(x_2-\langle \mathbf{x}\rangle)^2+(x_3-\langle \mathbf{x}\rangle)^2 = \textcolor{black}{2r/3}\right\rbrace
$$
in $X$ are stationary. In what follows, we focus exclusively on detailed-balance stationary states as they admit an explicit description. \textcolor{black}{However, we have found in a follow-up investigation that the results on detailed-balance stationary states in the following (sub)sections fully generalize to all stationary states \cite{Marc_arxiv}.}
\subsection{Stability conditions}\label{SS3A: stability conditions}
As mentioned earlier, the stationary states of a dynamical system do not always correspond to a robust equilibrium. To characterize the stability of a state, we study how a perturbed state $\mathbf{x}^\star+\boldsymbol{\epsilon}$ evolves and in particular, whether it converges back to $\mathbf{x}^\star$ or not. To this end, we assume the perturbation $\boldsymbol{\epsilon}$ to be sufficiently small such that the dynamics are determined by the linearized system around $\mathbf{x}^\star$ as
\begin{equation}\label{eq: linearized system}
\frac{d(\mathbf{x}^\star+\boldsymbol{\epsilon})}{dt} \approx J(\mathbf{x}^\star)\boldsymbol{\epsilon}\quad\text{with~}(J(\mathbf{x}^\star))_{ij} = \frac{d^2x_i}{dtdx_{j}}\Big\vert_{\mathbf{x}=\mathbf{x}^\star}
\end{equation}
where $J(\mathbf{x}^\star)$ is the Jacobian of system \eqref{eq: dynamical system} at $\mathbf{x}^\star$. If this Jacobian is (positive) negative definite, it implies directly that the stationary state $\mathbf{x}^\star$ is \textcolor{black}{\emph{linearly (un)stable}, which in turn implies (no) asymptotic stability. As we only consider stability criteria following from the linearized system \eqref{eq: jacobian elements}, we will further omit `linear' and simply write stable and unstable. If the} Jacobian is semi-definite instead, the linearization is not sufficient to determine the stability of $\mathbf{x}^\star$ and other techniques are required. 
\\
Restricting our analysis to the detailed-balance stationary states \eqref{eq: detailed-balance states}, we can further simplify the linearized system \eqref{eq: linearized system} and characterize certain stable stationary states. Here, we present a first stability result for system \eqref{eq: dynamical system}:
\begin{proposition}[full consensus/dissensus stability]\label{prop: full consensus/dissensus}
On any graph, the following states (if they exist) are stable stationary states of system \eqref{eq: dynamical system}
$$
\begin{cases}
(r<0):\textup{~the full consensus state~} (x^\star_i-x^\star_{j\sim i})=0\\
(r>0):\textup{~the full dissensus state~}(x^\star_i-x^\star_{j\sim i})=\pm\sqrt{r} 
\end{cases}
$$
For $r>0$ the full consensus state is unstable.
\end{proposition}
\textbf{Proof:~} See Section \ref{SS3B: laplacian form}.$~\square$\\
Proposition \ref{prop: full consensus/dissensus} only gives a rough picture of the stability of system \eqref{eq: dynamical system}, but it does illustrate clearly how the local dynamics are manifest in the global dynamics: the fact that dissensus is stable for each link ($1$-dimensional pitchfork system) locally while consensus is unstable, is observed globally as well. In the following section we refine this picture and show that the interconnected system also supports different types of stable states which are not simply inherited from the local dynamics. In particular, we find that for $r>0$ in the range between full consensus (and thus instability) and full dissensus (and thus stability) there may be \emph{stable mixed states} with both types of links present. As consensus links cannot exist stably for the local dynamics (see pitchfork dynamics, Section \ref{SS1A: three perspectives}), the existence of these stable mixed states is necessarily a feature of the system as a whole.
\subsection{Laplacian form of the linearized system}\label{SS3B: laplacian form}
Before continuing our analysis, we introduce some more information about the graph (network) on which the system takes place. By $G=(\mathcal{N},\mathcal{L})$ we will denote a graph with a set of $N$ nodes $\mathcal{N}$, and a set of $L$ links $\mathcal{L}\subseteq \mathcal{N}\times\mathcal{N}$ that connect pairs of distinct nodes, written as $i\sim j$ or $(i,j)\in\mathcal{L}$. We assume the graph to be finite and connected, i.e. with at least one path connecting each pair of nodes. Any graph $G$ has a corresponding $(N\times N)$ Laplacian matrix $Q$, with entries defined by
$$
(Q)_{ij} \triangleq \begin{cases} 
d_i &\text{~if~} i=j\\
-1  &\text{~if~} i\sim j\\
0   &\text{~otherwise},
\end{cases}
$$
where the degree $d_i$ of a node $i$ equals the number of neighbours of $i$ in $G$. The Laplacian matrix is just one among several matrix representations, but it is known to have close relations to many important graph properties \cite{Mohar, Merris, Chung} and appears in the formulation of diffusion processes on a given graph \cite{pvm_PA}.
\\
In the case of our system, the Laplacian matrix appears when calculating the Jacobian $J(\mathbf{x}^\star)$ of the system around some detailed-balance stationary state $\mathbf{x}^\star$. From their definition in \eqref{eq: linearized system}, we find that the entries of the Jacobian equal
\begin{equation}\label{eq: jacobian elements}
(J(\mathbf{x}^\star))_{ij} = \begin{cases}
d_ir - 3\sum_{j\sim i}(x^\star_i-x^\star_j)^2 &\text{~if~}i=j\\
-r + 3(x_i^\star-x_j^\star)^2 &\text{~if~}i \sim j\\
0 &\text{~otherwise}.
\end{cases}
\end{equation}
We let $\mathcal{L}_{=}\triangleq\lbrace i\sim j: x^\star_i-x^\star_j=0\rbrace$ and $\mathcal{L}_{\neq}\triangleq \lbrace i\sim j : x^\star_i-x^\star_j = \pm\sqrt{r} \rbrace$ denote the links of $G$ over which there is consensus, respectively dissensus in $\mathbf{x}^\star$. Correspondingly, we define the Laplacians $Q_{=}$ and $Q_{\neq}$ of the subgraphs of $G$ restricted to the consensus, respectively dissensus links \footnote{In other words, $Q_=$ and $Q_{\neq}$ are the Laplacian matrices of the graphs $G_{=}=(\mathcal{N},\mathcal{L}_{=})$ and $G_{\neq}=(\mathcal{N},\mathcal{L}_{\neq})$ containing only the consensus or dissensus links}; these matrices satisfy $Q=Q_{=}+Q_{\neq}$ since $\mathcal{L}=\mathcal{L}_{=}\cup\mathcal{L}_{\neq}$ holds. This subgraph decomposition allows the Jacobian to be expressed as follows:
\begin{lemma}\label{lemma: laplacian form of jacobian}
The Jacobian $J(\mathbf{x}^\star)$ of system \eqref{eq: dynamical system} at a detailed-balance solution $\mathbf{x}^\star$ with consensus and dissensus links $\mathcal{L}_{=}$ and $\mathcal{L}_{\neq}$ can be written as
\begin{align}\label{eq: sum of laplacians}
J(\mathbf{x}^\star) &= r(Q-3Q_{\neq}) \\
&= r(3Q_=-2Q)\nonumber\\
&= r(Q_{=}-2Q_{\neq})\nonumber
\end{align}
\end{lemma}
\textbf{Proof:~} Identity \eqref{eq: sum of laplacians} follows directly from the elementwise expression \eqref{eq: jacobian elements} and the definition of Laplacian matrices.$~\square$\\
Lemma \ref{lemma: laplacian form of jacobian} implies that the stability problem for detailed-balance stationary states comes down to characterizing  the spectrum of a difference of Laplacian matrices and, in particular, the positivity/negativity of its spectrum. 
\\
An important result about the Laplacian matrix of a connected graph is that it is positive semidefinite (i.e. non-negative eigenvalues) with a single zero eigenvalue corresponding to the constant eigenvector \cite{Mohar, Merris}. As the state space $X$ is orthogonal to the constant vector (by conservation of average), the Laplacian is thus effectively \emph{positive definite}. This observation leads to a direct proof of the stability result from Section \ref{SS3B: laplacian form}. \\
\textbf{Proof of Proposition \ref{prop: full consensus/dissensus}:~}If $\mathbf{x}^\star$ is the full consensus stationary state, we have that $\mathcal{L}_==\mathcal{L}$ and thus $J(\mathbf{x}^\star) = rQ$, which is positive definite if $r>0$ ($\mathbf{x}^\star$ is unstable) and negative definite when $r<0$ ($\mathbf{x}^\star$ is stable). If $\mathbf{x}^\star$ is the full dissensus state on the other hand, we find $\mathcal{L}_{\neq}=\mathcal{L}$ such that $J(\mathbf{x}^\star)=-2rQ$ which is negative definite for $r>0$ ($\mathbf{x}^\star$ is stable).$~\square$\\
While Proposition \ref{prop: full consensus/dissensus} is a direct result of the relation \eqref{eq: sum of laplacians} between the Jacobian and the Laplacian matrix of the graph on which the system takes place, the result does not depend on the specific structure of $G$ but only on the properties of the Laplacian matrix in general. The specific structure will play an important role in the case of mixed stationary states.
\subsection{Stability via effective resistances}\label{SS3C: effective resistances}
Somewhat surprisingly, the stability of mixed stationary states can be characterized in terms of the \emph{effective resistance}. The effective resistance was originally defined in the context of electrical circuit theory, but has found its way into graph theory through various applications such as random walks \cite{Doyle}, distance functions \cite{Klein}, graph embeddings \cite{Fiedler_book} and, more recently, graph sparsification \cite{Spielman}. The effective resistance $\omega_{ij}$ between a pair of nodes $i$ and $j$ in a graph $G$ can be defined as
\begin{equation}\label{eq: definition of effective resistance}
\omega_{ij} = (e_i-e_j)^TQ^\dagger(e_i-e_j),
\end{equation}
with $Q^\dagger$ the pseudoinverse of the Laplacian of $G$. For more intuition into the effective resistance, we refer the readers to \cite{Ghosh, Dorfler_electrical}, where expression \eqref{eq: definition of effective resistance} is derived starting from the electrical circuit equations. One of the important properties of the effective resistance is that it determines a \emph{metric} between the nodes of $G$ \cite{Klein}, where a small effective resistance between a pair of nodes indicates that these nodes are essentially close and `well connected', while a large effective resistance indicates the opposite. For instance, for a pair of linked nodes $i\sim j$ the extreme values for effective resistance correspond to $\omega_{ij}=2/N$ for the complete graph (i.e. very well connected) and $\omega_{ij}=1$ for a tree graph (i.e. poorly connected).
\\
We can now continue to characterize the stability of detailed-balance stationary states in the $r>0$ regime. From Proposition \ref{prop: full consensus/dissensus} we know that in full dissensus the system is stable while full consensus is unstable. Here, we provide a partial answer to the stability question for \emph{mixed detailed-balance states} with both consensus and dissensus links. In particular, we consider the case where a single consensus link is added to an otherwise full dissensus state; in this case, the stability depends on which link the consensus takes place:
\begin{theorem}[single consensus link stability]\label{th: single link stability} For system \eqref{eq: dynamical system} with $r>0$ on any graph $G$, the mixed stationary state $\mathbf{x}^\star$ with a single consensus link $\mathcal{L}_{=}=\lbrace (i,j)\rbrace$ satisfies
$$
\begin{cases}
\omega_{ij}<2/3:~\mathbf{x}^\star\textup{~is stable}\\
\omega_{ij}>2/3:~\mathbf{x}^\star\textup{~is unstable}
\end{cases}
$$
\end{theorem}
\textbf{Proof:} The proof is given in Appendix \ref{A1: single link stability} and is based on Lemma \ref{lemma: laplacian form of jacobian} and a new approach to bound the eigenvalues of a difference of Laplacian matrices.$~\square$\\
Theorem \ref{th: single link stability} states that a single consensus link state \emph{can} be stable, depending on the effective resistance of the consensus link $i\sim j$. Importantly, the criteria in Theorem \ref{th: single link stability} are \emph{tight} (except for a single point). If the effective resistance of the consensus link is high, i.e. if $i$ and $j$ are not well-connected, the state will not be stable. As mentioned before, an extreme example of the effective resistance is the case of tree graphs, where each pair of nodes has only a single link between them with no other possible paths such that $\omega_{ij}=1$. Generally, a large effective resistance indicates `bridge links', i.e. links between nodes which have few (or long) parallel paths between them (see example in Section \ref{SS6C: barbell graph}). Adding more parallel paths between $i$ and $j$ will gradually reduce the effective resistance until $\omega_{ij}=2/3$ is crossed, at which point the corresponding mixed state turns stable. In other words, \emph{bridge-like links} with few alternative paths in parallel cannot sustain consensus, while links with many alternative parallel paths can. 
\\
The answer to the initial question whether mixed stationary states can be stable is thus \emph{yes}, with the important condition that the consensus occurs between well-connected nodes. The proof of Theorem \ref{th: single link stability} is easily adapted to provide a condition for mixed states with several consensus links:
\begin{proposition}[mixed stationary state stability]\label{prop: mixed state stability} For system \eqref{eq: dynamical system} with $r>0$ on any graph $G$, the mixed stationary state $\mathbf{x}^\star$ with consensus links $\mathcal{L}_{=}$ satisfies
\begin{equation}\label{eq: mixed stability criteria}
\begin{cases}
\sum_{\mathcal{L}_{=}}\omega_{ij}<2/3:~\mathbf{x}^\star\textup{~is stable}\\
\max_{\mathcal{L}_{=}}\omega_{ij}>2/3:~\mathbf{x}^\star\textup{~is unstable}
\end{cases}
\end{equation}
\end{proposition}
\textbf{Proof:} See Appendix \ref{A2: mixed link stability}. \\
While Proposition \ref{prop: mixed state stability} is applicable to all mixed stationary states, the stability criteria are not tight like the criteria of Theorem \ref{th: single link stability}. Indeed, there are generally many detailed-balance states $\mathbf{x}^\star$ on a graph which satisfy neither of the criteria \eqref{eq: mixed stability criteria} and for which Proposition \ref{prop: mixed state stability} thus does not apply. \textcolor{black}{As discussed in Section \ref{SS6C: barbell graph}, one of the consequences of Theorem \ref{th: single link stability} and Proposition \ref{prop: mixed state stability} seems to be that in networks with a community structure, the stable states will generally contain more dissensus links \emph{between} different communities than \emph{within}. This would result in a higher similarity of node states within each of the communities, compared to an expected bias between the communities, which is an attractive modeling feature e.g. in the context of social cleavage \cite{Friedkin}. Crucially however, the results in Section \ref{SS6B: complete graph} show that within each of the communities, a certain level of dissensus is still expected to occur - the so-called spontaneous symmetry breaking described in \cite{Franci} - which can also be explained based on effective resistances in the graph, as shown in \cite{Marc_arxiv}.}
\\~\\
\emph{To summarize}, we studied the stationary states of system \eqref{eq: dynamical system} and identified the detailed-balance states \eqref{eq: detailed-balance states} as a subset of all possible stationary states. The characterization of the Jacobian matrix around detailed-balance stationary states as a difference of Laplacian matrices (Lemma \ref{lemma: laplacian form of jacobian}) enables a characterization of the stability in terms of the effective resistance. Most importantly, we find a \emph{tight stability condition} for states with a single consensus link (Theorem \ref{th: single link stability}) as well as more general, but less tight conditions for any mixed stationary state (Proposition \ref{prop: mixed state stability}). \textcolor{black}{In follow-up work \cite{Marc_arxiv}, we have found that all these results generalize to the setting of system \eqref{eq: odd non-linear system} with any odd coupling function $f$, and for all stationary states (using a suitable reformulation).}
%%%%%
%
%%%%%
\section{Exact solutions}\label{S4: exact solutions}
On certain networks, the stationary states $\mathbf{x}^\star$ of system \eqref{eq: dynamical system} can coincide with eigenvectors of the network Laplacian $Q$. As developed in detail in \cite{Prasse} for contagion dynamics, this allows for an \emph{exact solution} of the state evolution. Applied to our system, we find the following result:
\begin{theorem}[Exact solution]\label{th: exact solution}
If system \eqref{eq: dynamical system} on a graph $G$ has a stationary state $\mathbf{x}^\star\in X$ which is also a Laplacian eigenvector satisfying $Q\mathbf{x}^\star=\mu \mathbf{x}^\star$, then the exact solution for initial state $\mathbf{x}_0 = \alpha_0\mathbf{x}^\star$ and $r>0$ is given by
\begin{equation}\label{eq: exact solution}
\mathbf{x}(t,\mathbf{x}_0) = \mathbf{x}_0\left(\alpha_0^2 - \left(\alpha_0^2-1\right)e^{-2\mu r t}\right)^{-1/2}
\end{equation}
In particular, the system will reach the stationary state $\lim_{t\rightarrow \infty} \mathbf{x}(t,\mathbf{x}_0)=\mathbf{x}^\star$.
\end{theorem}
\textbf{Proof:~}see Appendix \ref{A3: exact solution}.$~\square$\\
In other words, Theorem \ref{th: exact solution} states that if the subspace $Z\subset X$ spanned by an eigenvector $\mathbf{z}$ of the Laplacian matrix contains a stationary state of system \eqref{eq: dynamical system}, then any initial condition in $Z$ allows for an exact solution \footnote{It follows directly from the proof in Appendix \ref{A3: exact solution} that if a stationary state $\mathbf{x}^\star$ is in a larger dimensional eigenspace of the Laplacian instead, i.e. $\mathbf{x}^\star=\sum c_k \mathbf{z}_k$ for some eigenvectors $\mathbf{z}_k$ of the Laplacian $Q$ with the same eigenvalue, that system \eqref{eq: dynamical system} can still be solved as in Theorem \ref{th: exact solution} for an initial state $\mathbf{x}_0=\alpha_0\mathbf{x}^\star$.}. Moreover, as $\mathbf{x}_0\in Z$ implies that $\mathbf{x}(t,\mathbf{x}_0)\in Z$, the subspace $Z$ is a \emph{positive invariant set} for the dynamics. For $r<0$, solution \eqref{eq: exact solution} still holds as long as $\vert\alpha_0\vert<1$.
\\
The question remains for which graphs there exist stationary states of system \eqref{eq: dynamical system} which are also Laplacian eigenvectors. In other words, we are looking for graphs for which there exists a state $\mathbf{x}^\star$ that satisfies
\begin{equation}\label{eq: definition of eigenstates}
\mu {x}^\star_i = \sum_{j\sim i}(x^\star_i-x^\star_j) \text{~~and~~} \sum_{j\sim i}r(x^\star_i-x^\star_j) = \sum_{j\sim i}(x^\star_i-x^\star_j)^3
\end{equation}
for all $i$. We will further refer to the states that satisfy \eqref{eq: definition of eigenstates} as \emph{eigenstates} of our system; regarding our system as a map $\phi_t : \mathbf{x}_0 \mapsto \mathbf{x}(t,\mathbf{x}_0)$, we find that $\phi_t(\mathbf{x}^\star) = \alpha_t \mathbf{x}^\star$ for these vectors, similar to the definition of eigenvectors for linear maps.\\
\emph{Example:} An elementary example of an eigenstate can be found for system \eqref{eq: dynamical system} on a pair of  connected nodes, i.e. $G=K_2$. The corresponding Laplacian matrix $\left(\begin{smallmatrix}1&-1\\-1&1\end{smallmatrix}\right)$ has a single non-constant eigenvector equal to $\mathbf{z}=(1,-1)^T$ with corresponding eigenvalue $\mu=2$. Scaling this eigenvector as $\mathbf{x}^\star = \sqrt{r}/2(1,-1)^T$ yields a detailed-balance stationary state, indicating that $\mathbf{x}^\star$ is an eigenstate of system \eqref{eq: dynamical system}. Consequently, the system can be solved exactly for $K_2$ consistent with the fact that we can solve the pitchfork normal form exactly, as shown in Figure \ref{fig: pitchfork}. In the following subsection we describe how simple examples like this two-node graph can be used as a starting point to construct new examples.
%%%
%
%%%
\subsection{Graphs with external equitable partitions}\label{SS4A: EEPs}
\begin{figure}[h!]
    \centering
    \includegraphics[scale=0.42]{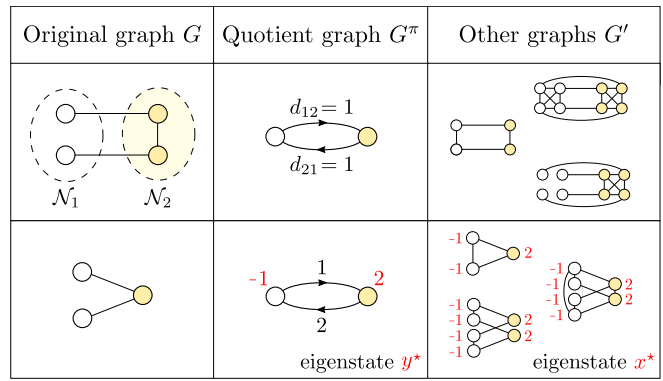}
    \caption{Illustration of external equitable partitions (EEPs), quotient graphs and the construction of eigenstates. In the first row, a partition $\pi$ of the nodes of $G$ (in two colours) is shown. Since each yellow node has one white neighbour and vice versa, this partition is an EEP and the corresponding quotient graph $G^{\pi}$ with directed, weighted links is given. In the third column, a number of other graphs $G'$ with EEPs are give for which $G^{\pi}$ is again the quotient graph. In the second row, another instance of $G,G^{\pi}$ and $G'$ is given, together with an example of how an eigenstate $\mathbf{y}^\star$ (satisfying \eqref{eq: definition generalized eigenstates}) on the quotient graph can be used to find eigenstates $\mathbf{x}^\star$ (satisfying \eqref{eq: definition of eigenstates}) on graphs $G'$ with EEPs.}
    \label{fig: equitable partitions}
\end{figure}

In the study of network dynamics and Laplacian matrices, an important type of graph symmetry are equitable partitions \cite{Schaub, Clery}. A partition $\pi$ of a graph divides the nodes of $G$ into $K$ disjoint groups $\mathcal{N}_1,\dots,\mathcal{N}_K\subseteq \mathcal{N}$ and is called an \emph{external equitable partition} (EEP) if all nodes in a group have the same number of links $d_{km}$ to all {other} groups, in other words 
\begin{equation}\label{eq: definition quotient graph weights}
\vert \lbrace v\sim i:v\in\mathcal{N}_m\rbrace\vert=\vert\lbrace v\sim j:v\in\mathcal{N}_m\rbrace\vert\triangleq d_{km}
\end{equation}
for all $i,j\in\mathcal{N}_{k\neq m}$. If $G$ has an external equitable partition $\pi$, its structure at the partition level can be summarized by the \emph{quotient graph} $G^{\pi}$. This quotient graph has node set $\lbrace 1,\dots,K\rbrace$ corresponding to the node groups of $G$ and a set of weighted, directed links $\overrightarrow{\mathcal{L}}$ that connect node group pairs $(k,m)$ between which there exist links in $G$, and with link weights $d_{km}$ for the link going from $k$ to $m$, and $d_{mk}$ for the link going from $m$ to $k$. Some examples of equitable partitions and quotient graphs are given in Figure \ref{fig: equitable partitions}. For more details on equitable partitions and their relation to dynamical systems, we refer the reader to \cite{Schaub, Clery}.
\\
The concept of external equitable partitions will allow us to construct eigenstates of system \eqref{eq: dynamical system} on graph $G$ based on eigenstates on its quotient graph $G^{\pi}$. Since $G^{\pi}$ is generally a directed and weighted graph, we generalize the definition of eigenstates to this setting as
\begin{align}\label{eq: definition generalized eigenstates}
 \sum_{m\sim k}d_{km}(y^\star_k-y^\star_m) &= \mu y^\star_k\quad\text{~and}\\ 
\sum_{m\sim k} d_{km}r(y^\star_k-y^\star_m) &= \sum_{m\sim k}d_{km}(y^\star_k-y^\star_m)^3.\nonumber
\end{align}
for all $k$. In Appendix \ref{A4: quotient graph} we show that if a vector $\mathbf{y}^\star$ satisfies \eqref{eq: definition generalized eigenstates} on $G^{\pi}$ then the corresponding vector $\mathbf{x}^\star$ with entries $x_i^\star=y_k^\star$ for $i\in\mathcal{N}_k$ will also satisfy \eqref{eq: definition of eigenstates} on $G$. As a result, we find that
\begin{proposition}\label{prop: eigenstates and EEPs}
For a graph $G$ with external equitable partition $\pi$, any eigenstate $\mathbf{y}^\star$ of system \eqref{eq: dynamical system} on $G^{\pi}$ has a corresponding eigenstate $\mathbf{x}^\star$ on $G$.
\end{proposition}
\textbf{Proof:~}See Appendix \ref{A4: quotient graph}.$~\square$\\
Proposition \ref{prop: eigenstates and EEPs} is a powerful tool for constructing examples of graphs with eigenstates. Indeed, starting from a (directed, weighted) graph $G$ with an eigenstate $\mathbf{y}^\star$ we can construct many examples of graphs $G'$ for which $G$ is a quotient graph, i.e. $G=G'^{\pi}$ with respect to an EEP $\pi$ of $G'$, and for which there thus exists an eigenstate $\mathbf{x}^\star$. In this construction, any node $k$ in $G$ can be replaced by a set $\mathcal{N}_k$ of nodes in $G'$ which can be interconnected in any desired way, and where the nodes from $\mathcal{N}_k$ are then given $d_{km}$ links to nodes in $\mathcal{N}_m$, which requires that the identity $\vert\mathcal{N}_k\vert d_{km} = \vert\mathcal{N}_m\vert d_{mk}$ holds for all pairs of partitions. This construction and the corresponding relation between eigenstates is illustrated in Figure \ref{fig: equitable partitions}.
\\~
\\
External equitable partitions arise from the notion \eqref{eq: definition quotient graph weights} of equivalence between nodes of a network, which is based on (local) symmetries between the neighbourhoods of the nodes. This relation between symmetry and dynamics closely resembles the perspective on network dynamics developed by Golubitsky, Stewart \emph{et al.} \cite{Golubitsky_groupoid_formalism, Stewart_Groupoids_and_patterns} and generalized in \cite{Golubitsky_Torok_multiple_arrows, Deville_Lerman, Nijholt_fibrations}. Their framework for general directed, labelled graphs and non-linear couplings focuses on \emph{local} symmetries (which carry the structure of a so-called `groupoid') of the graph, and studies how these symmetries are manifest in dynamics that respect the network structure. Our result that eigenstates $y^\star$ on the quotient graph $G^{\pi}$ can be `lifted' to eigenstates on the graph $G$ can be directly understood in context of this framework as a relation between (EEP) symmetry and dynamics, as shown in \cite{Aguiar}. These results are complementary to other works that study the effect of \emph{global} symmetries (i.e. graph automorphisms which carry the structure of a `group') on dynamical properties \cite{Ashwin_swift, Golubitsky_advances_in_symmetry}.
%%%
%
%
%%%
\section{Basins of attraction}\label{S5: basins of attraction}
Another classical question in (non-linear) dynamics is to determine which initial conditions lead to which stationary states. More specifically, the problem consists of characterizing the \emph{basins of attraction} $W(\mathbf{x}^\star)$ of the stationary states $\mathbf{x}^\star$, which are subsets of the state space defined as \footnote{An alternative formulation starts by defining the equivalence relation $\mathbf{v}\equiv \mathbf{w} \Leftrightarrow \lim_{t\rightarrow\infty}\mathbf{x}(t,\mathbf{v})=\lim_{t\rightarrow\infty}\mathbf{x}(t,\mathbf{w})$, which determines the basins of attraction as {equivalence classes} of the state space $X$; since $\lim_{t\rightarrow\infty}\mathbf{x}(t,\mathbf{x}^\star)=\mathbf{x}^\star$ trivially, these classes can be represented as $[\mathbf{x}^\star]$ and we have $W(\mathbf{x}^\star)=[\mathbf{x}^\star]$.}
$$
W(\mathbf{x}^\star) \triangleq \left\lbrace \mathbf{x}_0\in X: \lim_{t\rightarrow \infty}\mathbf{x}(t,\mathbf{x}_0)=\mathbf{x}^\star\right\rbrace,
$$
where we recall that $\mathbf{x}(t,\mathbf{x}_0)$ is the system state at time $t$ with initial condition $\mathbf{x}(0)=\mathbf{x}_0$. The basins of attraction are positive invariant sets, since for each $\mathbf{w}\in W(\mathbf{x}^\star)$ we have $\mathbf{x}(t,\mathbf{w})\in W(\mathbf{x}^\star)$ for all $t>0$, and determine a \emph{partition} of the state space of system \eqref{eq: dynamical system} as
\begin{equation}\label{eq: basins of attraction partition}
X = \bigcup_{\text{stationary~}\mathbf{x}^\star}W(\mathbf{x}^\star) 
\end{equation}
with any pair of distinct basins disjoint $W(\mathbf{x}^\star)\cap W(\mathbf{x}'^\star)=\emptyset$. Less formally, expression \eqref{eq: basins of attraction partition} captures the intuitive fact that for any initialization, system \eqref{eq: dynamical system} will converge to some stationary state. In general, it can be difficult to determine the basins of attraction for a non-linear system, but in the case of system \eqref{eq: dynamical system} we can use the additional properties of the dynamics (see e.g. Section \ref{SS1A: three perspectives}) to find some partial characterization. For instance, using the non-increasing property of the potential $V$, we know that the basin of attraction of a stationary point $\mathbf{x}^\star$ can only contain states of a higher potential, i.e. that $V(\mathbf{w})\geq V(\mathbf{x}^\star)$ for each $\mathbf{w}\in W(\mathbf{x}^\star)$.
\\
When system \eqref{eq: dynamical system} takes place on a \emph{tree graph $T$}, we can say even more about the basins of attraction. In Propositions \ref{prop: stationary=detailed-balance on trees} and \ref{prop: full-dissensus=only stable on trees} in Section \ref{SS6A: loopless networks}, we will show that all stationary states are detailed-balance stationary states and that among these states only the full dissensus states are stable. By \eqref{eq: potential in terms of dissensus links}, this means that all stable states on $T$ have the same minimal potential $V_{\min}=-r^2L/4$. Moreover, from the non-increasing property of the potential we find that there is a critical potential $V_c \triangleq -(L-1)r^2/4$ which determines a transition between states (with potential $V>V_c$) which in principle could be in the basin of attraction of any stationary state, and states (with potential $V<V_c$) which can only be in the basin of attraction of a stable state. We find the following characterization of the basins of attraction in the sub-critical regime:
\begin{proposition}[attraction basins on trees]\label{prop: attraction basins on trees}
For system \eqref{eq: dynamical system} on tree graphs, the state space region with potential lower than the critical potential $X\vert_{V<V_c}$ can be partitioned into basins of attraction of just the stable stationary states
\begin{equation}\label{eq: subcritical state space partition}
\left(X = \bigcup_{\text{stable~}\mathbf{x}^\star} W(\mathbf{x}^\star)\right)\Bigg\vert_{V<V_c}
\end{equation}
Furthermore, the basins of attraction in this region are given by
\begin{align}\label{eq: subcritical basin of attraction}
\left(W(\mathbf{x}^\star) = \Big\lbrace \mathbf{x}\in X : (x_i-x_j)(x^\star_i-x^\star_j)>0 ~\forall i\sim j\right\rbrace\Big)\Big\vert_{V<V_c}.
\end{align}
\end{proposition}
\textbf{Proof:~}See Appendix \ref{A5: basins of attraction}.$~\square$
%%%
%
%%%
\section{Examples and modeling observations}\label{S6: examples}
In the previous sections, we focused on the technical analysis of system \eqref{eq: dynamical system} and in particular on its stationary states. In the rest of the article, we study the system on a number of prototypical networks. We give a qualitative  description of the system solutions and suggest how certain properties might be useful when considering our system as a complex systems model.
%%%
%
%%%
\subsection{System on loopless networks}\label{SS6A: loopless networks}
On a loopless network, or \emph{tree graph} $T$, several of the earlier results are simplified or hold with less restrictions. Firstly, since the graph contains no loops, any assignment of $\lbrace 0,\pm\sqrt{r}\rbrace$ to the links of $T$ is possible; this amounts to $3^{N-1}$ possible detailed-balance stationary states (as any connected tree has $N-1$ links). Moreover, condition \eqref{eq: stationary states system} for stationarity implies \eqref{eq: detailed-balance states} for the case of tree graphs, and thus:
\begin{proposition}\label{prop: stationary=detailed-balance on trees}
For system \eqref{eq: dynamical system} on tree graphs, all stationary states (satisfying \eqref{eq: stationary states system}) are detailed-balance stationary states (satisfying \eqref{eq: detailed-balance states}).
\end{proposition}
\textbf{Proof:~}See Appendix \ref{A7: tree graphs}.$~\square$\\
Furthermore, the effective resistance between any pair of nodes of a tree graph is equal to the geodesic distance between these nodes \cite{Klein} which means that for linked nodes $i\sim j$ we have $\omega_{ij}=1$. Consequently, by Theorem \ref{th: single link stability} and Proposition \ref{prop: mixed state stability} we find that the stability of system \eqref{eq: dynamical system} is given by
\begin{proposition}\label{prop: full-dissensus=only stable on trees}
For system \eqref{eq: dynamical system} on tree graphs and $r>0$, the full dissensus state is stable while all other stationary states are unstable.
\end{proposition}
\textbf{Proof:} From Proposition \ref{prop: stationary=detailed-balance on trees}, the fact that $\omega_{ij}=1$ for all links $i\sim j$ in $T$ and Theorem \ref{th: single link stability}, it follows that any stationary state $\mathbf{x}^\star$ with a consensus link, i.e. with $\mathcal{L}_{=}$ non-empty, has $\max_{\mathcal{L}_{=}}{\omega_{ij}} = 1 >2/3\Rightarrow \mathbf{x}^\star$ is unstable. The stability of the full dissensus state follows from Proposition \ref{prop: full consensus/dissensus}.$~\square$\\
One consequence of Proposition \ref{prop: full-dissensus=only stable on trees} is that the proportion of stable stationary states on a tree equals $(2/3)^{N-1}$ which vanishes exponentially fast for larger trees.\\
As discussed in the previous section in Proposition \ref{prop: attraction basins on trees}, we also have some information about the basins of attraction for tree graphs.
%%%
%
%%%
\subsection{Balanced opinion formation in the complete graph}\label{SS6B: complete graph}
In the complete graph $K_N$, every node is connected to all $(N-1)$ other nodes, making it the densest possible graph. Moreover, it means that the graph contains a high level of symmetry in the sense that no two nodes are distinguishable from their connections to other nodes, which greatly simplifies the description of detailed-balance stationary states. Since any three nodes in $K_N$ form a triangle, the only stationary states $(x_i^\star,x_j^\star,x_k^\star)$ these nodes can achieve is some permutation of
$$
(0,0,0) \text{~or~}(0,0,\pm\sqrt{r}).
$$
As a consequence, the stationary states at the level of the full graph must be some permutation of the state
$$
(\underbrace{0,\dots,0}_{\#(N-V)},\underbrace{\sqrt{r},\dots,\sqrt{r}}_{\# V})
$$
with $V$ entries equal to $\sqrt{r}$ and $N-V$ entries equal to $0$. In other words, all stationary states are parametrized by the number $V\in[0,N]$ of $\sqrt{r}$-states; since there are ${{N}\choose{V}}$ ways to choose $V$ such nodes, the complete graph has $2^N$ detailed-balance stationary states (by the binomial theorem). Moreover, if we use the degree of freedom provided by the average state to fix the state of an (arbitrary) reference node to $x^\star_i=0$, the (rescaled) state parameter $v=V/N$ will be related to the average state value by $\langle \mathbf{x}^\star \rangle = \pm v \sqrt{r}$. For the stability of the stationary states, we find the following result:
\begin{proposition}\label{prop: stable states on K_N}
For system \eqref{eq: dynamical system} on the complete graph $K_N$ with $r>0$, the detailed-balance stationary states $\mathbf{x}^\star$ satisfy
$$
\begin{cases}
v\in(1/3,2/3) : \mathbf{x}^\star\text{~is stable}\\
v\notin[1/3,2/3] : \mathbf{x}^\star\text{~is unstable} 
\end{cases}
$$
\end{proposition}
\textbf{Proof:~}See Appendix \ref{A7: tree graphs}.$~\square$\\
This characterization of the (stable) stationary states on the complete graph is very interesting from a modeling point of view. First, in any detailed-balance stationary state the nodes of $G$ are split into two groups with $V$ and $(N-V)$ nodes, respectively. Within each of these groups the nodes are in consensus, while between the groups there is dissensus. Furthermore, Proposition \ref{prop: stable states on K_N} states that the size of the two groups needs to be \emph{balanced} in stable states, i.e. the group sizes can differ by at most $N/3$ and neither of the groups can dominate the full graph. Figure \ref{fig: complete graph} below illustrates the findings of Proposition \ref{prop: stable states on K_N} in the bifurcation diagram of $K_N$. 
\\
In contrast to loopless graphs, Proposition \ref{prop: stable states on K_N} shows that a non-vanishing proportion of $2/3$ of all detailed-balance stationary states are stable on the complete graph.
\begin{figure}[h!]
    \centering
    \includegraphics[scale=0.7]{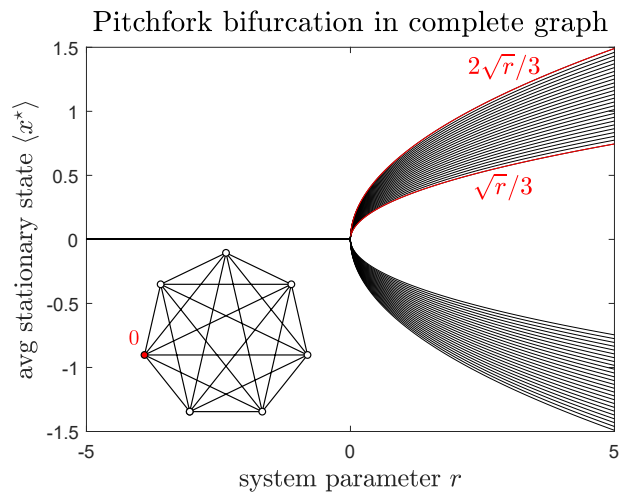}
    \caption{Bifurcation diagram of system \eqref{eq: dynamical system} on the complete graph $K_N$ with $N=75$ nodes. For $r<0$ the consensus state is the only stable stationary state. For $r>0$, the stable detailed-balance stationary states are parametrized by $v\in(1/3,2/3)$ with corresponding average state value $\langle \mathbf{x}^\star\rangle = \pm v\sqrt{r}$ when fixing an arbitrary node to $x^\star_i=0$. The diagram resembles the bifurcation diagram (Figure \ref{fig: pitchfork}) of the pitchfork bifurcation normal form.}
    \label{fig: complete graph}
\end{figure}
\\~
\\
Assuming the framework of \emph{opinion dynamics} \cite{Castellano, Franci, Bizyaeva}, where nodes play the role of individuals in a population with states $x_i(t)$ recording their preference in the range between a certain opinion A with $x_i=+\sqrt{r}/2$, or an opposing opinion B if $x_i=-\sqrt{r}/2$  (rivaling political party, competing product, etc.), we might interpret these results as follows: for $r<0$ any difference between initial individual preferences will be disappear from the network until the population reaches a global consensus where all individuals agree. This qualitative behaviour is studied in various contexts like engineering \cite{Saber} or social sciences \cite{Castellano, Proskurnikov}, and can also be reproduced by the simpler diffusion dynamics $dx_i/dt=\sum_{j\sim i}(x_j-x_i)$. For $r>0$ on the other hand, an atypical stationary distribution emerges in system \eqref{eq: dynamical system} where instead of reaching global consensus, the population splits into two groups adhering to different opinions. Moreover, the stability condition $v\in (1/3,2/3)$ guarantees that neither of these groups can be too dominant in the population, i.e. that there is a balanced coexistence of opinions. This qualitative behaviour is observed in real social systems, where it is often called social cleavage or polarization \cite{Friedkin}. \textcolor{black}{As noted in \cite{Franci} it is remarkable that a fully interconnected system with indistinguishable nodes (agents) can exhibit spontaneous symmetry breaking into a state with distinct groups of nodes. The analysis in \cite{Franci} however explains how this behaviour is expected for a broad class of dynamical systems}
\\ Finally, we remark that the complete graph should be seen as a prototype for more general `dense graphs', and that our qualitative description should hold approximately for dense random graphs like, for instance, Erd\H{o}s-R\'{e}nyi random graphs with high link probability $p$, as a result of concentration of measure \cite{Bandeira}. Importantly, the above description of the equilibrium behaviour of system \eqref{eq: dynamical system} on $K_N$ does not take any non detailed-balance states into account, for which we might observe very different types of stable states.
\\~\\
\textcolor{black}{An equivalent system consisting of pitchfork bifurcation normal forms on the complete graph was analysed by Aronson \emph{et al.} in \cite{Aronson_Josephson_junctions} as a model of coupled Josephson junctions. In contrast to our ad-hoc derivation, they make a principled \emph{equivariant analysis} of the system dynamics and derive the two-group stationary states from the fact that this fully interconnected system has $S_n\times \mathbbm{Z}_2$ symmetry (permuting nodes $\times$ sign change). The same stability conditions as Proposition \ref{prop: stable states on K_N} are noted in \cite{Aronson_Josephson_junctions} based on calculations of the system Jacobian for the pitchfork bifurcation normal form (similar to our proof), however there is no suggestion as to how this stability result might generalize to less symmetrical systems. In a sense, this broader view on the relation between (in)stability and structure is exactly what Lemma \ref{lemma: laplacian form of jacobian} and Theorem \ref{th: single link stability} in the present work and some (stronger) results in \cite{Marc_arxiv} build towards. This illustrates the complementarity between the approach and tools of \cite{Franci, Bizyaeva} (symmetric graphs but general systems) and our approach (general graphs but specific system).}
%%%%
%
%%%%
\subsection{Biased opinion formation in the barbell graph}\label{SS6C: barbell graph}
\begin{figure}[h!]
    \centering
    \includegraphics[scale=0.7]{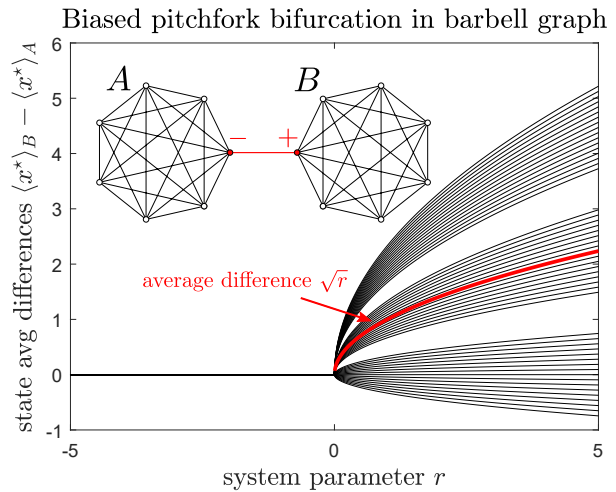}
    \caption{Difference in average stationary state value of the two complete subgraphs for system \eqref{eq: dynamical system} on the barbell graph $B_{2N}$ with $2N=48$ nodes. For $r<0$ the consensus state is the only stable stationary state. For $r>0$, the bridge link $i\sim j$ needs to be a dissensus link for stability (by Theorem \ref{th: single link stability}), which allows to fix $x^\star_i=-\sqrt{r}/2$ and $x^\star_j=\sqrt{r}/2$. The stable stationary states are then parametrized by $v_A,v_B\in(1/3,2/3)$ for each of the complete components independently, and the \emph{forced bias} over the bridge link results in a non-zero average difference of $\sqrt{r}$ between the state averages.}
    \label{fig: barbell graph}
\end{figure}

The barbell graph $B_{2N}$, illustrated in Figure \ref{fig: barbell graph}, consists of two complete graphs joined by a single link $i\sim j$. Similar to the complete graph, the high number of symmetries in each of the individual complete parts allows for a compact description of the stationary states. In fact, the detailed-balance stationary states on $B_{2N}$ can be parametrized as the stationary state on two `independent' complete graphs, i.e. with $V_A,V_B\in[0,N]$ denoting the number of nodes with a different value from $x^\star_i,x^\star_j$ in the two complete graphs respectively. The average state value in the complete graphs is then related to their respective (scaled) parameters as $\langle \mathbf{x}^{\star}\rangle_A = x^\star_i \pm v_A\sqrt{r}$ and similarly for $\langle \mathbf{x}^\star\rangle_B=x^\star_j\pm v_B\sqrt{r}$. Comparing the average state value between the two components yields
$$
\langle \mathbf{x}^{\star}\rangle_B-\langle \mathbf{x}^{\star}\rangle_A = (x^\star_i-x^\star_j) \pm (v_B\pm v_A)\sqrt{r}.
$$
When restricting to stable stationary states, we find that the bridge link in the barbell graph has effective resistance $\omega_{ij}=1$ which implies that this bridge link must have dissensus in all stable states. Furthermore, by Proposition \ref{prop: stable states on K_N} we have that the stationary states in the complete graphs are stable for $v_A,v_B\in(1/3,2/3)$. Assuming that $x^\star_j>x^\star_i$ we then find that the difference between the average \emph{stable} state values in both components equals
$$
\langle \mathbf{x}^\star\rangle_B - \langle \mathbf{x}^\star\rangle_A = \left[1\pm (v_B\pm v_A)\right]\sqrt{r} \text{~with~}v_A,v_B\in(1/3,2/3).
$$
In Figure \ref{fig: barbell graph} this finding is illustrated in the bifurcation diagram of $B_{2N}$, which clearly shows the non-zero difference that exists between both complete components for $r>0$.
\\~
\\
The stable detailed-balance states on the barbell graph are again interesting from a modeling perspective. In the setting of consensus dynamics, we again find that for $r<0$ all individuals converge to a common opinion. For $r>0$ on the other hand, a balanced coexistence of opinions will be established \emph{within} each of the complete graphs separately but, importantly, with a non-zero bias existing on average \emph{between} the components. In other words, within the dense subgraphs the opinions coexist without either opinion dominating the other, while a difference will exist between the subgraphs. This qualitative behaviour might seem interesting if we think of the barbell graph as a prototypical example of a graph consisting of dense groups of nodes which are sparsely interconnected in between, a structure commonly known as \emph{assortative communities}. In this setting, one might expect such opinion biases to exist between communities rather than within due to a different level of coordination or communication, and a model similar to our system might help explain the underlying mechanisms. 
\\~\\
\textcolor{black}{
We remark that the stationary states in this example could also be derived based on the $(S_{N-1}\wr \mathbbm{Z}_2)\times \mathbbm{Z}_2$ symmetry of the system (permuting non-bridge nodes within a complete component and/or interchanging components $\times$ sign change), where `$\wr$' denotes the wreath product between groups \cite{Wells}. From the theory of equivariant dynamics \cite{Golubitsky_advances_in_symmetry, Franci}, we know that certain stationary states will be associated to subgroups of these system symmetries. This analysis would not provide any stability information however.}
%%%
%
%%%

\section{Related result: synchronization of phase oscillators}\label{S7: synchronization}
A famous example of non-linear dynamics on networks are systems of interacting phase oscillators. The underlying idea is that many natural (herds/shoals of animals, groups of neurons, etc.) and man-made (power grids, electrical oscillators, etc.) systems can be modeled effectively as a population of oscillators which establish some form of synchronization due to interactions \cite{Strogatz_chaos, Arenas}. The periodic behaviour of a single entity is abstracted as an oscillator whose state $\theta_i(t)\in\mathbb{R}/2\pi$ cycles periodically according to a natural frequency $f_i$. These oscillators are then interconnected in a network $G$, with a coupling function $h$ driving the phases of adjacent entities to a common value as
\begin{equation}\label{eq: system of coupled oscillators}
\frac{d\theta_i}{dt} = f_i + \mu\sum_{j\sim i}h(\theta_i-\theta_j)
\end{equation}
with the coupling strength $\mu$ as a system parameter. The easiest example of a periodic, odd coupling function is the sine function. Similar to how our non-linear system \eqref{eq: dynamical system} is the $3^{\text{rd}}$ order Taylor approximation for any odd coupling function $f$ (on $\mathbb{R}$), the sine function can be seen as the $1^{\text{st}}$ order Fourier expansion for any periodic odd coupling function $h$ (on $\mathbb{R}/2\pi$). System \eqref{eq: system of coupled oscillators} with $h(x)=\sin(x)$ on the complete graph is also known as the \emph{Kuramoto model} and is widely studied in the context of synchronization, see for instance \cite{Strogatz_review, Arenas, Dorfler_synch_review}. 
\\
One of the key features that motivates the study of interacting oscillator systems is that the oscillators exhibit \emph{synchronization} for certain parameter ranges of $\mu$ and $\lbrace f_i\rbrace$ on certain graph structures. The onset of various types of synchronization (phase, frequency, chimera, etc.) has been studied extensively in these systems, and is used as a theoretical explanation for observed synchrony in many real-world systems. Here, we mention a specific result about coupled oscillators which is similar to our main result Theorem \ref{th: single link stability}.
\\
A particular notion of synchronization is ``frequency synchronization with $\gamma$-cohesive phases", which is defined as the (rotating) state $\boldsymbol{\theta}^\star$ where all oscillators rotate at the same instantaneous frequency $d\theta_i^\star/dt=f^\star$ and where the phase differences between adjacent oscillators in the network satisfy $(\theta^\star_i-\theta^\star_j)\operatorname{mod} 2\pi\leq \gamma$ for all $i\sim j$. We will call such a state \emph{$\gamma$-synchronized}. In \cite{Dorfler_smartgrids} the authors propose to study for which choices of natural frequencies $\mathbf{f}=(f_1,\dots,f_N)^T$ this type of synchrony can occur. Their interesting finding is that for many graphs (certain extremal graphs, and dense sets of random graphs) the following criterion
$$
\max_{i\sim j}\left\vert (e_i-e_j)Q^\dagger \mathbf{f}\right\vert\leq \sin(\gamma)
$$
is a sufficient condition for system \eqref{eq: system of coupled oscillators} on a graph with Laplacian $Q$ and sinusoidal coupling, to have a $\gamma$-synchonized state. In particular, this implies the known result that that system \eqref{eq: system of coupled oscillators} with a constant natural frequency $\mathbf{f}=\alpha \mathbf{1}$ can have a $\gamma$-synchronized state, for any $\gamma$. Moreover, if the natural frequencies are equal for all but one pair of connected nodes $i$ and $j$, which differ by $\vert f_i-f_j\vert=c$, then the synchrony criterion becomes $c\leq \sin(\gamma)/\omega_{ij}$, i.e. the difference $c$ is upper-bounded by the inverse of the effective resistance. In other words, starting from the {constant} frequency distribution for which there is {synchrony} possible, and changing a {single link} to have a frequency difference of $c$, then synchrony is conserved \emph{depending on the effective resistance} of the respective link. More specifically, a small (large) effective resistance will admit a large (small) phase difference. While the setting of \cite{Dorfler_smartgrids} is very different, this result is reminiscent of Theorem \ref{th: single link stability}, and a further investigation of this similarity might be worthwhile.
%%%
%
%
%%%
\section{Conclusion}\label{S8: conclusion}
In this paper, we have introduced and studied a  non-linear dynamical system on networks  inspired by the pitchfork bifurcation. Our analysis is motivated by different interpretations of the system as a collection of interdependent pitchfork systems, a gradient dynamical system for a potential composed of interacting double-well potentials and finally as the dominating behaviour for more general non-linear systems with odd coupling functions. In a certain sense system \eqref{eq: dynamical system} is the `simplest' dynamical system of a broad class of non-linear systems (e.g. with general odd coupling functions $f$ in \eqref{eq: odd non-linear system}, or general symmetric potentials $V$). The choice to study equations \eqref{eq: dynamical system} specifically is thus the outcome of a wish to implement a model with more complexity than simple linear models, while wielding Occam's razor.
\\
Our technical analysis mainly focused on the equilibrium behaviour of the system. The bifurcation from a single stationary state to a myriad of possible stationary states and in particular their stability, provides a clear picture of how the simple local dynamical rule in our system gives rise to interesting global phenomena. Specifically, as a main technical result (Theorem \ref{th: single link stability}) we found stability conditions that depend on the full structure of the network, as captured by their dependence on the \emph{effective resistance}. Our further analysis of the system includes the identification of exact solutions for certain graphs, which include graphs with external equitable partitions, and the description of basins of attraction for loopless networks. 
\\
Finally, we looked at the system on a number of prototypical graphs and describe some interesting qualitative behaviour of the solutions. On the complete and barbell graph, our results suggest an interpretation of the system as an opinion dynamic model: in one parameter regime the system is driven to a global consensus state, while the stable states in the other regime are characterized by a balanced bipartition of opinions (states) in dense components and with an overall non-zero bias between sparsely connected components, that grows as the system goes deeper in the parameter regime. These results support system \eqref{eq: dynamical system} as a rich model for complex systems allowing to identify unexpected bridges between network properties, like the effective resistance, and dynamical ones, which could trigger future advances in the more general study of non-linear systems on networks.
%%%%%%%%%%%%%%%%%%
% ACKNOWLEDGEMENTS
%%%%%%%%%%%%%%%%%%
\section*{Acknowledgements}
The authors would like to thank Christian Bick, Alain Goriely and Michael Schaub for helpful discussions, suggestions and references, Bastian Prasse for suggesting to study the exact solutions as in \cite{Prasse}, and Marc Homs Dones for carefully reading the manuscript. K.D. was supported by The Alan Turing Institute under the EPSRC grant EP/N510129/1. R.L. acknowledges support from the Flagship European Research Area Network (FLAG-ERA) Joint Transnational Call “FuturICT 2.0”

%%%%%%%%%
%APPENDIX
%%%%%%%%%
\appendix
\section*{Appendix}

\section{Proof of single consensus link stability}\label{A1: single link stability}
We prove Theorem \ref{th: single link stability} about the stability of mixed detailed-balance stationary states $\mathbf{x}^\star$ with a single consensus link $\mathcal{L}_{=}=\lbrace i,j\rbrace$. \\
\textbf{Proof:} By Lemma \ref{lemma: laplacian form of jacobian}, we know that the system Jacobian at $\mathbf{x}^\star$ equals
$$
J(\mathbf{x}^\star) = 3r\left((\mathbf{e}_i-\mathbf{e}_j)(\mathbf{e}_i-\mathbf{e}_j)^T - 2/3Q\right)
$$
with unit vectors $(\mathbf{e}_i)_k=1$ iff $k=i$ and 
where the first term equals $Q_{=}$ which is the (rank-one) Laplacian of a single-link graph. We introduce the abbreviations $J(\mathbf{x}^\star)=J$ and $\lambda_{\max}=\lambda_{\max}(J)$. In what follows, we show that the sign of $\lambda_{\max}$ equals the sign of $(\omega_{ij}-2/3)$ which by the linearization of system \eqref{eq: dynamical system} then implies Theorem \ref{th: single link stability}.
\\
By the \emph{Courant Minmax principle} for real symmetric matrices, we know that
\begin{equation}\label{eqA: courant minmax}
\lambda_{\max} \geq \mathbf{x}^TJ\mathbf{x} \text{~for all~} \mathbf{x}\in\mathbb{R}^N\text{~with~}\Vert \mathbf{x}\Vert = 1.
\end{equation}
In particular, we thus know that for $\tilde{\mathbf{x}}=Q^\dagger(\mathbf{e}_i-\mathbf{e}_j)$ and $\mathbf{x}=\tilde{\mathbf{x}}/\Vert \tilde{\mathbf{x}}\Vert$ we have
\begin{equation}
\lambda_{\max} \geq 3r\mathbf{x}^T\left((\mathbf{e}_i-\mathbf{e}_j)(\mathbf{e}_i-\mathbf{e}_j)^T - 2/3Q\right)\mathbf{x}
\end{equation}
Taking definition \eqref{eq: definition of effective resistance} of the effective resistance into account, this can be rewritten as
\begin{equation}\label{eq: lambda lower bound}
\lambda_{\max} \geq 3r\omega_{ij}/\Vert \tilde{\mathbf{x}}\Vert^2 (\omega_{ij}-2/3).
\end{equation}
Since the factor $3r\omega_{ij}/\Vert \tilde{\mathbf{x}}\Vert^2$ is positive for $r>0$ we have that $(\omega_{ij}-2/3)\geq0 \Rightarrow \lambda_{\max} \geq 0$.
\\
To upper-bound the largest eigenvalue, we start again from the Courant Minmax principle \eqref{eqA: courant minmax} and use the fact that equality is attained for some vector $\mathbf{y}$, i.e. such that
$$
\exists \mathbf{y}\text{~with~}\Vert \mathbf{y}\Vert = 1\text{, such that~} \lambda_{\max}=\mathbf{y}^TJ\mathbf{y}
$$
Writing out the full expression for the Jacobian, we find
$$
\lambda_{\max} = 3r \left(\left[(\mathbf{e}_i-\mathbf{e}_j)^T\mathbf{y}\right]^2 - 2/3\mathbf{y}^TQ\mathbf{y}\right)
$$
If we denote the quadratic product by $P=\mathbf{y}^TQ\mathbf{y}$, then we know that $\mathbf{y}$ satisfies $\Vert \mathbf{y}\Vert =1, \mathbf{y}^TJ\mathbf{y} = \lambda_{\max}$ and $\mathbf{y}^TQ\mathbf{y} = P$. If we \emph{relax} the first two conditions on $\mathbf{y}$, we get an upper-bound of the form
\begin{equation}\label{eqA: optimization in terms of z}
    \lambda_{\max} \leq 3r\max_{\mathbf{z}: \mathbf{z}^TQ\mathbf{z} = P}\left\lbrace  \left[ (\mathbf{e}_i-\mathbf{e}_j)^T\mathbf{z}\right]^2-2/3\mathbf{z}^TQ\mathbf{z} \right\rbrace.
\end{equation}
This is a valid upper-bound since $\mathbf{y}$ is in the domain of $\mathbf{z}$, but the domain for $\mathbf{z}$ is larger than $\mathbf{y}$ alone and can thus potentially achieve a larger objective value. Next, we introduce the decomposition of the (positive semidefinite) Laplacian as $Q=S^TS$ and its pseudoinverse as $Q^{\dagger} = S^{\dagger T}S^{\dagger}$, where the $(N-1)\times N$ matrices $S,S^\dagger$ have entries $(S)_{ij}=(z_j)_i\sqrt{\mu_j}$ and $(S^\dagger)_{ij} = (z_j)_i\sqrt{1/\mu_j}$ and thus satisfy $S^T S^{\dagger}=I-\mathbf{u}\mathbf{u}^T/N$ and $S^{\dagger}S^T=I_{N-1}$ \footnote{For an interpretation of this matrix $S$ as the vertex vector matrix of a simplex corresponding to the graph $G$, see \cite{Fiedler_book, krl_simplex}}, with all-one vector $\mathbf{u}$. With the change of variable $\tilde{\mathbf{z}} = S\mathbf{z} \Leftrightarrow \mathbf{z}=S^{\dagger T}\tilde{\mathbf{z}}$ we can then rewrite \eqref{eqA: optimization in terms of z} as
$$
\lambda_{\max} \leq 3r\max_{\tilde{\mathbf{z}}:\Vert \tilde{\mathbf{z}}\Vert^2=P}\left\lbrace \left[\tilde{\mathbf{z}}^T S^{\dagger}(\mathbf{e}_i-\mathbf{e}_j)\right]^2 -2/3\Vert \tilde{\mathbf{z}}\Vert^2 \right\rbrace,
$$
which is solved by taking $\tilde{\mathbf{z}}$ parallel to $S^{\dagger }(\mathbf{e}_i-\mathbf{e}_j)$, i.e. as $\tilde{\mathbf{z}} = S^{\dagger }(\mathbf{e}_i-\mathbf{e}_j)\sqrt{P/\omega_{ij}}$. Introducing the definition of the effective resistance in the form $\omega_{ij}=\left[ S^\dagger(\mathbf{e}_i-\mathbf{e}_j)\right]^2$ we then find \begin{equation}\label{eq: lambda upper bound}
\lambda_{\max} \leq 3rP(\omega_{ij}-2/3).
\end{equation}
Since the factor $3rP$ is positive for $r>0$, we have that $(\omega_{ij}-2/3)\leq 0\Rightarrow \lambda_{\max}\leq 0$. Combining inequalities \eqref{eq: lambda lower bound} and \eqref{eq: lambda upper bound}, we indeed find that $\operatorname{sign}(\lambda_{\max})=\operatorname{sign}(\omega_{ij}-2/3)$ which implies Theorem \ref{th: single link stability}.$~\square$
%%%%
%
%
%%%%
\section{Proof of Proposition \ref{prop: mixed state stability}}\label{A2: mixed link stability}
As in Appendix \ref{A1: single link stability}, we will bound the eigenvalues of the Jacobian matrix $J(\mathbf{x}^\star)$ starting from Lemma \ref{lemma: laplacian form of jacobian}. For consensus links $\mathcal{L}_{=}$ we have
$$
J(\mathbf{x}^\star) = 3r\sum_{(i,j)\in\mathcal{L}_{=}}(\mathbf{e}_i-\mathbf{e}_j)(\mathbf{e}_i-\mathbf{e}_j)^T - 2rQ.
$$
By Courants Minmax principle we can write
\begin{align}
\lambda_{\max} &\geq c(\mathbf{e}_a-\mathbf{e}_b)^TQ^\dagger JQ^\dagger(\mathbf{e}_a-\mathbf{e}_b)\nonumber\\
&= 3rc\left(\omega_{ab}^2 - 2/3\omega_{ab} \right) \nonumber
\\
&+ \sum_{(i,j)\in\mathcal{L}\backslash\lbrace(a,b)\rbrace}c\left((\mathbf{e}_a-\mathbf{e}_b)^TQ^\dagger(\mathbf{e}_i-\mathbf{e}_j)\right)^2\nonumber
\\
&\geq 3rc\omega_{ab}(\omega_{ab}-2/3)\label{eq: mixed lambda lower bound}
\end{align}
for any $(a,b)\in\mathcal{L}_{=}$ and with positive normalizing constant $c=\Vert Q^\dagger(\mathbf{e}_a-\mathbf{e}_b)\Vert^{-2}$. Since \eqref{eq: mixed lambda lower bound} holds for any consensus link, it is clear that $\max_{\mathcal{L}_{=}}\omega_{ab}>2/3\Rightarrow\lambda_{\max}>0$ proving the first stability condition.
\\
For the upper-bound, we proceed similar to Appendix \ref{A1: single link stability} and let $\mathbf{y}$ be the normalised vector such that $\mathbf{y}^TJ\mathbf{y}=\lambda_{\max}$ and define $P=\mathbf{y}^TQ\mathbf{y}$. We can then write
\begin{align}
\lambda_{\max} &= \mathbf{y}^TJ\mathbf{y} = 3r\mathbf{y}^TQ_{=}\mathbf{y} - 2r\mathbf{y}^TQ\mathbf{y}\nonumber
\\
&\leq 3r \max_{\mathbf{z}:\mathbf{z}^TQ\mathbf{z}=P}\left\lbrace \sum_{(i,j)\in\mathcal{L}_{=}}\left[(\mathbf{e}_i-\mathbf{e}_j)^T\mathbf{z}\right]^2 - 2/3P\right\rbrace\nonumber
\\
&\leq 3r\sum_{(i,j)\in\mathcal{L}_{=}}\left\lbrace\max_{\mathbf{z}:\mathbf{z}^TQ\mathbf{z}=P}\left[(\mathbf{e}_i-\mathbf{e}_j)^T\mathbf{z}\right]^2 - 2/3P\right\rbrace\nonumber
\\
&= 3rP\left(\sum_{(i,j)\in\mathcal{L}_{=}}\omega_{ij}-2/3\right)\label{eq: mixed lambda upper bound}
\end{align}
where the maximum in the third line is achieved by $\mathbf{z}=S^{\dagger}(\mathbf{e}_i-\mathbf{e}_j)(\sqrt{P/\omega_{ij}})$ for each consensus link $i\sim j$. The upper-bound \eqref{eq: mixed lambda upper bound} shows that if $\sum_{\mathcal{L}_{=}}\omega_{ij}<2/3\Rightarrow \lambda_{\max}<0$ proving the second stability condition.$~\square$
%%%
%
%
%%%

\section{Proof of exact solution}\label{A3: exact solution}
We prove Theorem \ref{th: exact solution} about exact solutions of system \eqref{eq: dynamical system} for certain initial conditions. In line with the approach of \cite{Prasse}, we will show that if at some time $t$, the state is of the form $\mathbf{x}(t)=\alpha(t) \mathbf{x}^\star$ with $\mathbf{x}^\star$ a stationary state of system \eqref{eq: dynamical system} parallel to an eigenvector of the Laplacian matrix with eigenvalue $\mu$ (i.e. $\mathbf{x}^\star$ is an eigenstate), the dynamic equations \eqref{eq: dynamical system} simplify to an equation for $\alpha(t)$. This proves that the subspace spanned by the vector (eigenstate) $\mathbf{x}^\star$ is a positive invariant subspace, as for any initial state $\mathbf{x}_0=\alpha_0\mathbf{x}^\star$ the solution is of the form $\mathbf{x}(t,\mathbf{x}_0)=\alpha(t)\mathbf{x}^\star$. Secondly, we show that the time-dependent coefficient $\alpha(t)$ can be solved exactly as the solution of a $1$-dimensional Bernoulli differential equation. \\
\textbf{Proof:~} The system equations for $\mathbf{x}(t)$ are given by
$$
\frac{dx_i}{dt} = \sum_{j\sim i}r(x_i-x_j) - (x_i-x_j)^3 \text{~for all~}i.
$$
If at some time $t$ the state is of the form $\mathbf{x}(t)=\alpha(t)\mathbf{x}^\star$, where $\mathbf{x}^\star$ is an eigenstate of the system satisfying conditions \eqref{eq: definition of eigenstates}, these equations become
\begin{equation}\label{eqA: system for alpha(t)x}
\frac{d[\alpha(t)x^\star_i]}{dt} = \sum_{j\sim i}r\alpha(t)(x^\star_i-x^\star_j) - \alpha^3(t)(x^\star_i-x^\star_j)^3.
\end{equation}
Since $\mathbf{x}^\star$ is a stationary state of the dynamics, we have that $\sum_{j\sim i}(x_i^\star-x_j^\star)^3=r\sum_{j\sim i}(x^\star_i-x^\star_j)$, and \eqref{eqA: system for alpha(t)x} simplifies to 
$$
\frac{d\alpha}{dt}x^\star_i = r\left(\alpha(t)-\alpha^3(t)\right)\sum_{j\sim i}(x^\star_i-x^\star_j).
$$
Next, as $\mathbf{x}^\star$ is an eigenvector of the Laplacian $Q$, i.e. with $\mu x^\star_i = \sum_{j\sim i}(x^\star_i-x^\star_j)$ we can rewrite this as
$$
\frac{d\alpha}{dt}x^\star_i = r\mu\left(\alpha(t)-\alpha^3(t)\right)x^\star_i.
$$
This shows that, for all $t'>t$ the solution will be of the form $\mathbf{x}(t',\mathbf{x}(t)) = \alpha(t')\mathbf{x}^\star$ and thus that $\left\lbrace c\mathbf{x}^\star \text{~for some~} c\in\mathbb{R}\right\rbrace$ is a positive invariant set for the dynamics \footnote{At this point, we have also implicitly assumed that $\alpha(t)$ is differentiable at $t$, i.e. that the state does not just cross the eigenspace determined by $x^\star$. This assumption is verified since solution \eqref{eqA: solution for alpha(t)} is consistent with system \eqref{eq: dynamical system}}.
\\
Furthermore, the equation for $\alpha(t)$ is a $1$-dimensional Bernoulli differential equation
$$
\frac{d\alpha}{dt} = r \mu (\alpha(t)-\alpha^3(t)),
$$
which can be solved by introducing $\beta=\alpha^{-2}-1$ with $d\beta/dt = -2\alpha^{-3}d\alpha/dt$. This yields a linear differential equation $\frac{d\beta}{dt} = -2r\mu\beta$ with solution $\beta(t)=\beta_0e^{-2r\mu t}$. Introducing the initial condition $\beta_0=\alpha_0^{-2}-1$ and changing the variable back to $\alpha(t)$ we find the solution
\begin{align}\label{eqA: solution for alpha(t)}
\alpha(t)&= \left[1+\left(\frac{1}{\alpha_0^2}-1\right)e^{-2\mu rt}\right]^{-1/2}\nonumber\\
&=\alpha_0\left(\alpha_0^2 - \left(\alpha_0^2-1\right)e^{-2\mu rt}\right)^{-1/2}.
\end{align}
which proves Theorem \ref{th: exact solution}. $\medskip\square$
\\
Note that this solution method only works when $(\beta(t)+1)^{-1/2}$ is well defined. When $r>0$ or $r<0$ and $\vert\alpha_0\vert<1$ this is satisfied; when $r<0$ and $\vert\alpha_0\vert>1$ on the other hand, we find that $t\geq\tfrac{1}{2\mu\vert{r}\vert}\ln\left(\frac{\alpha_0^2}{\alpha_0^2-1}\right)\Rightarrow \beta(t)\leq-1$ in which case solution \eqref{eqA: solution for alpha(t)} will not apply.
%%%%
%
%
%%%%
\section{Proof of Proposition \ref{prop: eigenstates and EEPs}}\label{A4: quotient graph}
Let $\mathbf{y}^\star$ be an eigenstate of the quotient graph $G^\pi$, and $\mathbf{x}$ a $(N\times 1)$ vector with entries $x_i=y^\star_k$ when $i\in\mathcal{N}_k$, defined on the nodes of $G$. For any $k$ and $i\in\mathcal{N}_k$ we can then write
\begin{align*}
\sum_{j\sim i}(x_i-x_j) &= \sum_{m=1}^K\sum_{\substack{j\sim i\\j\in\mathcal{N}_m}}(x_i-x_j)\\
&= \sum_{m\neq k}\sum_{\substack{j\sim i\\j\in\mathcal{N}_m}}(x_i-x_j)\\
&= \sum_{m\neq k}(y^\star_k-y^\star_m)\sum_{\substack{j\sim i\\j\in\mathcal{N}_m}}1\\
&= \sum_{m\sim k}d_{km}(y^\star_k-y^\star_m) = \mu y^\star_k = \mu x_i
\end{align*}
where we used definition \eqref{eq: definition quotient graph weights} of the link weight $d_{km}$. This illustrates that $\mathbf{x}$ is an eigenvector of $Q$ with eigenvalue $\mu$. Similarly, for any $k$ and $i\in\mathcal{N}_k$ we can write
\begin{align*}
\sum_{j\sim i}(x_i-x_j)^3 &= \sum_{m\neq k}\sum_{\substack{j\sim i\\j\in\mathcal{N}_m}}(x_i-x_j)^3\\
&= \sum_{m\neq k}(y_k^\star-y_m^\star)^3\sum_{\substack{j\sim i\\j\in\mathcal{N}_m}}1\\
&= \sum_{m\sim k}d_{km}(y^\star_k-y^\star_m)^3
\\
&= \sum_{m\sim k}rd_{km}(y^\star_k-y^\star_m)=\sum_{j\sim i}r(x_i-x_j)
\end{align*}
which shows that $\mathbf{x}$ is also a stationary state of system \eqref{eq: dynamical system}. Together, these identities show that the vector $\mathbf{x}$ constructed from the eigenstate $\mathbf{y}^\star$ of $G^\pi$ is an eigenstate of $G$. $~\square$
%%%%
%
%
%%%%
\section{Proof of attraction basins on trees}\label{A5: basins of attraction}
We prove Proposition \ref{prop: attraction basins on trees}, which characterizes the basins of attraction for system \eqref{eq: dynamical system} in the region of state space with potential below the critical potential $V_c=-(L-1)r^2/4$.\\
\textbf{Proof:~} The non-increasing property of the potential $\dot{V}\leq 0$ implies that any state $\mathbf{x}$ with $V(\mathbf{x})<V_c$ can only evolve to stationary states $\mathbf{x}^\star$ with $V(\mathbf{x}^\star)\leq V(\mathbf{x})<V_c$. Since this is only leaves the stable stationary states (i.e. with $\ell(\mathbf{x}^\star)=L$ and $V(\mathbf{x}^\star)=-Lr^2/4$), we have that any $\mathbf{x}$ with $V(\mathbf{x})<V_c$ can only converge to a stable stationary state. Furthermore, by partition \eqref{eq: basins of attraction partition} of the \emph{full} state space, every state is in the basin of attraction of exactly one stationary state, which thus means that $\mathbf{x}$ is in the basin of attraction of exactly one stable stationary state. Since this holds for every state with sub-critical potential, we get decomposition \eqref{eq: subcritical state space partition} of the \emph{sub-critical region} of state space into basins of attraction of the stable states.
\\
Next, we derive the specific form of the basins of attraction in the sub-critical region. We start by showing that for states $\mathbf{w}$ with $V(\mathbf{w})<V_c$ we have $\operatorname{sign}(x_i(t,\mathbf{w})-x_j(t,\mathbf{w}) = \operatorname{sign}(w_i-w_j)$ for all links $i\sim j$ and $t>0$. In other words, the state difference over a link can not `flip' its sign. Assume for contradiction that there is a link $i\sim j$ for which the sign is flipped. Since the state differences satisfy the differential equation \eqref{eq: system with link variables}, $(x_i-x_j)(t)$ must be continuous and thus we have that 
\begin{align*}
&\text{if~}(x_i-x_j)(t_0)>0\text{~and~}(x_i-x_j)(t_1)<0,\\
&\text{then~} \exists \tau \in [t_0,t_1]\text{~s.t.~}(x_i-x_j)(\tau) = 0
\end{align*}
The potential of the state $\mathbf{x}=\mathbf{x}(\tau)$ at this time can then be lower-bounded by
\begin{align*}
V(\mathbf{x}) &= \frac{1}{4}\sum_{\substack{a\sim b \\ \neq i\sim j}}(x_a-x_b)^2((x_a-x_b)^2-2r) 
\\
&+ \frac{1}{4}\underbrace{(x_i-x_j)^2((x_i-x_j)^2-2r)}_{=0}
\\
&\geq \sum_{\substack{a\sim b \\ \neq i\sim j}}\frac{-r^2}{4} = -\frac{(L-1)r^2}{4} = V_c
\end{align*}
where we used the fact that $z^2(z^2-2r)\geq -r^2$. This bound contradicts the fact that $V(\mathbf{x}(t,\mathbf{w}))\leq V(\mathbf{w})<V_c$ which means our assumption must be false. Thus, since the stable stationary states have $(x^\star_i-x^\star_j)=\pm\sqrt{r}$ and we know that $\operatorname{sign}(x_i(t,\mathbf{w})-x_j(t,\mathbf{w})) = \operatorname{sign}(w_i-w_j)$ for all $t>0$ and for $t\rightarrow \infty$ in particular, we know that the link-difference signs of a state in the sub-critical region of $X$ determines its corresponding stationary state. This proves \eqref{eq: subcritical basin of attraction} in Proposition \ref{prop: attraction basins on trees}. $~\square$
%%%%
%
%
%%%%
\section{Stationary states on tree graphs}\label{A7: tree graphs}
To prove Proposition \ref{prop: stationary=detailed-balance on trees}, we will make use of the fact that each (non-singleton) tree graph $T$ has at least one \emph{leaf node} $l$ (a node with degree $d_l=1$), and that removing a leaf node from $T$ results in another tree $T'=T\backslash \lbrace l\rbrace$.
\\
\textbf{Proof of Proposition \ref{prop: stationary=detailed-balance on trees}:} For $N=2$, the only tree $T$ consists of a pair of connected nodes $i\sim j$. System \eqref{eq: stationary states system} which determines the stationary states reduces to a single equation $r(x^\star_i-x^\star_j)=(x_i^\star-x^\star_j)^3$ which solves to a detailed-balance stationary state \eqref{eq: detailed-balance states}. We will prove the general case of Proposition \ref{prop: stationary=detailed-balance on trees} by induction on the number of nodes $N$, with base case $N=2$. Assuming the induction hypothesis holds for $N=K$, i.e. for all trees on $K$ nodes the stationary states are detailed-balance stationary states, we will now show that it holds for $N=K+1$ as well. 
\\
Let $T$ be a tree graph on $N=K+1$ nodes and $l$ one of its leaf nodes connected to one other node $n$. The stationary states $\mathbf{x}^\star$ of system \eqref{eq: dynamical system} on $T$ are found from $d\mathbf{x}^\star/dt=0$ which, for the leaf node, yields the equation 
\begin{equation}\label{eqA: equation for leaf node}
\frac{dx_l^\star}{dt} = 0 \Leftrightarrow
r(x^\star_l-x^\star_n) - (x^\star_l-x^\star_n)^3=0,
\end{equation}
which is only satisfied if $(x^\star_l-x^\star_n)\in\lbrace 0,\pm\sqrt{r}\rbrace$. In other words, for any (not necessarily detailed-balance) stationary state $\mathbf{x}^\star$ of system \eqref{eq: dynamical system} on $T$, the link difference $(x^\star_l-x^\star_n)$ must be either a consensus link or a dissensus link. The stationary state values on the other nodes $i\neq l$ are determined by the equations
\begin{equation}\label{eq: system on non-leaf nodes}
\begin{dcases}
\frac{dx^\star_i}{dt}=0 &\Leftrightarrow \sum_{j\sim i}r(x_i^\star-x_j^\star) + (x_i^\star-x_j^\star)^3 = 0\text{~for~}i\neq n\\
\frac{dx^\star_n}{dt}=0 &\Leftrightarrow \sum_{j\sim n, j\neq l}r(x_n^\star-x_j^\star) + (x_n^\star-x_j^\star)^3  \\
&\hphantom{\Leftrightarrow}+ r(x_n^\star-x_l^\star) - (x_n^\star-x_l^\star)^3= 0
\end{dcases}
\end{equation}
Since leaf node $l$ has degree $d_l=1$, the only equation where $x^\star_l$ appears is the balance equation for $x^\star_n$. Moreover, introducing $(x_n^\star-x^\star_l)\in\lbrace 0,\pm\sqrt{r}\rbrace$ in equation \eqref{eq: system on non-leaf nodes} for $x^\star_n$ eliminates $x^\star_l$ from the equations for $\mathbf{x}^\star$ altogether. From this, it follows that the stationary state $\mathbf{x}^\star$ can be determined from the stationary states ${\mathbf{x}^{\star}}'$ of the tree graph $T'=T\backslash\lbrace l\rbrace$ as $x^\star_i={x^{\star}}'_i$ for all $i\neq l$ (as they obey the same equations) and with $(x^\star_l-x^\star_n)\in\lbrace 0,\pm\sqrt{r}\rbrace$ (by the solution of equation \eqref{eqA: equation for leaf node}). By the induction hypothesis, ${\mathbf{x}^{\star}}'$ is the stationary state of system \eqref{eq: dynamical system} on a $K$-node tree graph and is thus a detailed-balance state. From this follows that $\mathbf{x}^\star$ will be a detailed-balance stationary state as well. $\medskip\square$ 
%%%
%
%
%%%
\section{Stability on the complete graph}\label{A6: complete graph}
We prove Proposition \ref{prop: stable states on K_N} which states that stationary states $\mathbf{x}^\star$ with $V\in(N/3,2N/3)$ are stable by explicitly calculating the spectrum of the Jacobian $J(\mathbf{x}^\star)$.
\\
\textbf{Proof:~} The stationary state $\mathbf{x}^\star$ partitions the set of nodes $\mathcal{N}$ of $K_N$ into two disjoint sets $\mathcal{V}$ with $\vert\mathcal{V}\vert=V$ and $\overline{\mathcal{V}}=\mathcal{N}\backslash\mathcal{V}$ such that for all $i\in\mathcal{V},j\in\overline{\mathcal{V}}$ we have $(x^\star_i-x_j^\star) = \pm\sqrt{r}$. In other words, all consensus links go between nodes within a same set, while dissensus links go between nodes of a different set. If we order the nodes as $\mathcal{V}={1,\dots,V}$ and $\overline{\mathcal{V}}={V+1,\dots,N}$ then the Laplacian matrices can be written as
\begin{align*}
&Q_{=} = \begin{pmatrix} VP_V & 0 \\ 0 & (N-V)P_{N-V} \end{pmatrix}\\
&Q_{\neq} = \begin{pmatrix} (N-V)I_V & -\mathbf{u}_V\mathbf{u}_{N-V}^T\\ -\mathbf{u}_{N-V}\mathbf{u}_V^T & VI_{N-V} \end{pmatrix}
\end{align*}
where $I_{*}$ and $\mathbf{u}_{*}$ denote the identity matrix and all-one vector of dimensions indicated by $*$, and with $P_V = I_V-\mathbf{u}_V\mathbf{u}_V^T/V$ the projector on the space orthogonal to $\mathbf{u}_V$ and similarly for $N-V$. The Jacobian of system \eqref{eq: dynamical system} at $\mathbf{x}^\star$ can then be calculated as $J = r(Q_{=}-2Q_{\neq})$. We will show that $J$ has four types of eigenvectors and, correspondingly, four types of eigenvalues.\\
\underline{Type 1:} Any vector of the form $\mathbf{z} = (\mathbf{z}_V, \mathbf{0}_{N-V})^T$ with $\mathbf{z}_V$ a $V$-dimensional vector that satisfies $\mathbf{z}_V^T\mathbf{u}_V=0$ and $\Vert \mathbf{z}_V\Vert = 1$, will give 
\begin{align*}
J\mathbf{z} &= r\left(Q_{=}\mathbf{z} - 2Q_{\neq}\mathbf{z}\right)\\
&= r(V - 2(N-V))\mathbf{z}\\
&= 3r(V-2/3N)\mathbf{z}
\end{align*}
which shows that $\mathbf{z}$ is an eigenvector of $J$ with eigenvalue $3r(V-2/3N)$. Since we can choose a basis of $V-1$ vectors orthogonal to $\mathbf{u}_V$ which all are of the form of $\mathbf{z}$, the Jacobian will have $V-1$ eigenvalues equal to $3r(V-2/3N)$.\\ \underline{Type 2:} The same approach based on vectors of the form $\mathbf{z}=(\mathbf{0}_{V},\mathbf{z}_{N-V})^T$ with $\mathbf{z}_{N-V}^T\mathbf{u}_{N-V}=0$ gives $(N-V-1)$ eigenvalues of $J$ equal to $-3r(V-N/3)$. \\
\underline{Type 3:} The third type of vector is given by $\mathbf{z}=(-\mathbf{u}_V/V,\mathbf{u}_{N-V}/(N-V))^T$ and has corresponding eigenvalue $0$, as $J\mathbf{z}=-2rN$.\\
\underline{Type 4:} Finally, the fourth type of vector is simply the constant vector $\mathbf{u}=(\mathbf{u}_V,\mathbf{u}_{N-V})$ with eigenvalue $0$, as $J\mathbf{u}=0$.
\\
By construction, we furthermore have that these $N$ vectors of the form $(\mathbf{z}_V,\mathbf{0}_{N-V}),(\mathbf{0}_V,\mathbf{z}_{N-V})$, $(-\mathbf{u}_{V}/V,\mathbf{u}_{N-V}/(N-V))$,$(\mathbf{u}_{V},\mathbf{u}_{N-V})$ are orthogonal and and thus determine eigenvectors of $J$. The eigenvalues of the Jacobian $J(\mathbf{x}^\star)$ are thus equal to
\begin{align*}
&\lbrace 0\rbrace^{(1)}, \lbrace-2rN\rbrace^{(1)},\left\lbrace -3r(V-N/3)\right\rbrace^{(N-V-1)},
\\
&\lbrace 3r(V-2N/3)\rbrace^{(V-1)}
\end{align*}
where superscripts denote the multiplicity of the eigenvalues. The first zero eigenvalue corresponds to the constant eigenvector $\mathbf{u}$ which is orthogonal to the state space $X$ and thus does not influence the stability of state $\mathbf{x}^\star$. The second zero eigenvalue is always negative for $r>0$. Finally, if $V\in(N/3,2N/3)$ then all other eigenvalues are negative and if $V\notin[N/3,2N/3]$ they are positive, which proves the (in)stability of the stationary states $\mathbf{x}^\star$ with respect to $V$. When $V\in\lbrace N/3,2N/3\rbrace$, the other eigenvalues become zero and the linearization method does not provide the necessary information to determine the stability of the corresponding states.$~\square$
%%%
%
%%%
%%%%%%%%%%%%%
% REFERENCES
%%%%%%%%%%%%%
\bibliography{bibliography_v2.bib}

\end{document}